\documentclass{eajam}
\usepackage{charter}
\usepackage[charter]{mathdesign}
\usepackage{bm}              
\def\bx{{\bm x}}
\def\by{{\bm y}}
\def\bz{{\bm z}}

\def\bg{{\bm g}}

\def\Rc{{\mathrm R}}

\def\Rc{{\mathrm R}}
\def\bP{{\mathrm  P}}
\def\bE{{\mathrm E}}

\def\bT#1{\left\|#1\right\|}
\def\BT#1{\left\|#1\right\|^2}

\def\bF #1{\| #1 \|_F}   
\def\BF #1{\| #1 \|_F^2}  

\allowdisplaybreaks[1]  

\usepackage{graphicx}
\usepackage{subfigure}
\usepackage{epstopdf} 

\usepackage{threeparttable} 

\usepackage{enumerate}  

\usepackage{algorithmicx,algorithm, algpseudocode}

\usepackage{color}
\begin{document}

\markboth{N.-C. Wu and C. Liu}{Randomized progressive iterative approximation}
\title{Randomized progressive iterative approximation for B-spline curve and surface fittings}

\author[N.-C. Wu and C. Liu]{Nian-Ci Wu \affil{1} and Chengzhi  Liu \affil{2}\comma\corrauth}
\address{\affilnum{1}\ School of Mathematics and Statistics, South-Central Minzu University, Wuhan 430074, China.\\
\affilnum{2}\ School of Mathematics and Finance, Hunan University of Humanities, Science and Technology, Loudi 417000,  China.}
%
%
\email{{\tt it-rocket@163.com} (C. Liu)}
%
\begin{abstract}
For large-scale data fitting, the least-squares progressive iterative approximation is a widely used method in many applied domains because of its intuitive geometric meaning and efficiency. In this work, we present a randomized progressive iterative approximation (RPIA) for the B-spline curve and surface fittings. In each iteration, RPIA locally adjusts the control points according to a random criterion of index selections. The difference for each control point is computed concerning the randomized block coordinate descent method. From geometric and algebraic aspects, the illustrations of RPIA are provided. We prove that RPIA constructs a series of fitting curves (resp., surfaces), whose limit curve (resp., surface) can converge in expectation to the least-squares fitting result of the given data points. Numerical experiments are given to confirm our results and show the benefits of RPIA.
\end{abstract}



\maketitle


 \section{Introduction}
Data fitting is a relevant problem in many applied domains, including computer aided design, computer graphics, data visualization, and many other fields. Progressive iterative approximation (PIA), also known as the geometric iteration method, is a class of
typical data fitting algorithms with clear geometric meaning, which avoids solving the whole linear systems directly and has been widely used in academic research and engineering practices. We refer to the  survey \cite{18LMD} for more details.

The PIA technique was respectively discovered by Qi et al. in 1975 and de Boor in 1979 and reignited by Lin et al. \cite{05LBW} in 2005.
PIA is generally categorized into interpolatory and approximate types.  There are several commonly employed interpolatory PIA, such as the
original PIA \cite{05LBW}, local PIA \cite{10Lin} and  weighted PIA \cite{10Lu,16ZGT}, which need the number of control points to equal that of data points.  Delgado and Pe\~{n}a \cite{07DP}  compared  the PIA convergence rates with different normalized totally positive (NTP) bases and showed that the normalized B-spline basis possesses the fastest convergence rate. In \cite{09MCK}, Martin et al. devised an equivalent PIA format with uniform periodic cubic B-spline. Chen et al. \cite{11CW} extended  the PIA property of univariate NTP basis to that of bivariate Bernstein basis over a triangle domain. In the function space formed by a NTP basis, PIA curve or surface fitting with a normalized B-spline basis has the fastest convergent speed \cite{07DP}.

In many applications, we are given more data points than what can be interpolated by a polynomial curve (resp., surface). In such cases, an approximating curve (resp., surface) will be needed. Such a curve (resp., surface) does not pass through the data points exactly; rather, it passes near them and captures the shape inherent to the data points.  The technique is known as the least-squares approximation, such as extended PIA \cite{11LZ}.

In a celebrated paper \cite{14DL}, Deng and Lin provided a least-squares PIA (LSPIA) with B-spline basis, which inherits the advantages of the classical PIA method and becomes quite useful in the shape modeling community \cite{19EL, 20HW, 18LLGZHS, 22RJ, 21Wang, 16ZGT}. Rios and J\"{u}ttler further excavated the algebraic property of LSPIA and proved that it is equivalent to a gradient descent method \cite{22RJ}. This approach generalizes broadly, for example,
to a linear dependent non-tensor product bivariate B-spline basis, which leads to a lower order surface fitting result \cite{18LLGZHS},
to a generalized B-spline basis \cite{16ZGT},
and to variational composite iterations \cite{19EL,20HW}.
 With the utilization of  the Schultz method, an LSPIA-type variant was developed by Ebrahimi and Loghmani in \cite{19EL}. We refer to it as SLSPIA for the rest of this paper. Very recently, Huang and Wang gave a PIA with memory for least-squares (MLSPIA) fitting as another improvement of LSPIA, where the information in the previous PIA step is necessary \cite{20HW}. 
LSPIA is recovered if one takes  some specific parameters in MLSPIA. As well as the update rule of the control points, many other works on extension for PIA are reported due to its promising performance and elegant mathematical property; see, for example \cite{13CDP,  20LHL, 07MMN, 12ONLKM, 06SW, 21Wang} and the references therein.


The approximate PIA variants, such as LSPIA \cite{14DL}, SLSPIA \cite{19EL}, and MLSPIA \cite{20HW},  are  global and need to adjust all of the control points simultaneously at each iteration.  Though local PIA \cite{10Lin} allows that partial control  points are updated, it fails to deal with the least-squares problem. Moreover, in the case of surface fitting,  all these methods need to operate the Kronecker product of two collocation matrices, which have a considerably large order.  In general, solving this way  is very demanding in terms of computational efficiency. To alleviate this issue, we give a new {\it approximate local PIA} method, which constructs a series of fitting curves or surfaces by adjusting  partial control points according to a random criterion and has a least-squares fitting result  to the given data points.

There are two major contributions of this work detailed as follows.

\begin{enumerate}[(1)]
\setlength{\itemindent}{0cm}
\vskip 0.5ex
\item {\sf The local least-squares solver.} Our method  is an approximate PIA algorithm. It allows the number of data points to be larger than that of control points and obtains a least-squares result in the limit sense. At each iteration step, our method only locally updates the control points concerning an index set and keeps the remaining control points unchanged.
\item{\sf Computing complexity.} Different from the traditional PIA surface fitting methods, our method does not need to operate the Kronecker product of two collocation matrices. It is proved equivalent to the solution of linear matrix equations from algebraic aspects.
\end{enumerate}

\vskip 0.5ex
We organize the remaining part of this paper as follows.
We first briefly review the LSPIA method for curve and surface fittings in Section \ref{Sec:2} and then introduce a randomized PIA  (this method is abbreviated as RPIA in the following) method  for curve  and surface  fittings in Section \ref{Sec:3}.
Next, Section \ref{Sec:4} shows its convergence.
Afterward, in Section \ref{Sec:5} some numerical examples are provided to demonstrate the theoretical results.
Finally, we end this paper with some conclusions in Section \ref{Con}.

\section{The LSPIA method}\label{Sec:2}
Rather than the traditional fitting methods directly based on the solution of a linear system, LSPIA generates a series of curves  (resp., surfaces)  to approximate the fitting curve (resp., surface) by a fixed parameter.

\subsection{The case of curves} Given an ordered point sequence $\left\{{\bm q}_j\right\}_{j=0}^m \subseteq \mathbb{R}^d$ $(d=2,3)$ to be fitted, a blending basis sequence $\left\{\mu_i(x):i\in [n]\right\}$  defined on $[0,1]$, and a real increasing sequence $\left\{x_j \right\}_{j=0}^m \subseteq [0,1]$, where the set $[\ell]:=\{0,1,\cdots,\ell\}$ for any positive integer $\ell$. For $k=0,1,2\cdots$, supposing that we have gotten the $k$th curve \begin{align*}
  \mathrm{C}^{(k)}(x) =\sum_{i=0}^n \mu_i(x) {\bm p}_i^{(k)},~x\in [x_0,~x_m]
\end{align*}
in which LSPIA iteratively approximates a target curve by updating the control point ${\bm p}_i^{(k)}$ for $i\in [n]$.
The difference, between ${\bm q}_j$ and the corresponding point on $\mathrm{C}^{(k)}(x)$, is defined by
\begin{align*}
  {\bm r}_{j}^{(k)} = {\bm q}_j - \mathrm{C}^{(k)}(x_j),~j\in[m],
\end{align*}
Taking a weighted sum of the differences, the $i$th adjusting vector is computed by
\begin{align*}
{\bm \delta}_i^{(k)} = \widetilde{\mu}  \sum_{j=0}^m \mu_i(x_j){\bm r}_{j}^{(k)}
\end{align*}
for $i\in[n]$ with $\widetilde{\mu}$ being a constant. Then the next curve is generated by
\begin{align*}
  \mathrm{C}^{(k+1)}(x)=\sum_{i=0}^n \mu_i(x) {\bm p}_i^{(k+1)}
\quad {\rm with} \quad
  {\bm p}_i^{(k+1)} = {\bm p}_i^{(k)} +  {\bm \delta}_i^{(k)}.
\end{align*}

Several well-known  LSPIA variants for curve fitting are obtained by setting appropriate choice for the adjusting vector.
\vspace{.5ex}
\begin{enumerate}[\quad\sffamily(1)]
       \setlength{\itemindent}{0cm}
       \item {\sf SLSPIA} \cite{19EL}. Based on the Schulz iterative method, i.e.,
        ${\bm Z}^{(k+1)} = 2{\bm Z}^{(k)}  - {\bm Z}^{(k)}{\bm A}{\bm Z}^{(k)},~ k=0,1,2,\cdots$,
        where ${\bm Z}^{(0)}=\widehat{\mu}{\bm A}^T {\bm A} {\bm A}^T$, ${\bm A}$ is the corresponding collocation matrix, and $\widehat{\mu}$ is a constant, for any $i\in [n]$, the $i$th adjusting vector is updated by
        \begin{align*} {\bm \delta}_{i}^{(k)} = \sum_{j=0}^m {\bm Z}_{ij}^{(k)} {\bm r}_{j}^{(k)}. \end{align*}
       \item {\sf MLSPIA} \cite{20HW}. This method needs to store and use the information of previous differences. For $i\in[n]$ and $k\geq 1$, by introducing three real weights $\omega$, $\gamma$, and $v$, the $i$th adjusting vector is computed by
          \begin{align*}
          {\bm \delta}_{i}^{(k)} = (1-\omega){\bm \delta}_{i}^{(k-1)}+\gamma\bar{\bm{\delta}}_{i}^{(k)}+(\omega-\gamma)\bar{\bm{\delta}}_{i}^{(k-1)},
          \end{align*}
          where $\bar{\bm{\delta}}_{i}^{(k)} = v \sum_{j=0}^m  \mu_i(t_j)  {\bm r}_{j}^{(k)}$. In particular,  one special case of MLSPIA reduces to  LSPIA when $\omega=\gamma=1$ and $v = \widetilde{\mu} $.
\end{enumerate}
\vspace{.5ex}
The list is by no means comprehensive and merely serves the purpose to illustrate the elasticity  of LSPIA, and here we omit an exhaustive review of the related pieces of literature \cite{18LCZ, 18LMD, 11LZ, 18LLGZHS, 16ZGT} and the references therein.

\subsection{The case of surfaces}Given an ordered point sequence $\left\{{\bm Q}_{hl}\right\}_{h=0,l=0}^{m,p} \subseteq \mathbb{R}^3$  to be fitted, and a real increasing sequence $\left\{y_l \right\}_{l=0}^p\subseteq [0,1]$. For $k=0,1,2,\cdots$, supposing that we have gotten the $k$th surface
\begin{align*}
  \mathrm{S}^{(k)}(x,y) =\sum_{i=0}^n \sum_{j=0}^n \mu_i(x) \mu_j(y) {\bm P}_{ij}^{(k)},~x\in [x_0,~x_m],~y\in [y_0,~y_p],
\end{align*}
where ${\bm P}_{ij}^{(k)}$ is the $(i,j)$th control point for $i,j\in[n]$. The $(h,l)$th difference is defined by
\begin{align*}
  {\bm R}_{hl}^{(k)} = {\bm Q}_{hl} - \mathrm{S}^{(k)}(x_h, y_l)
\end{align*}
for $h\in[m]$ and $l\in[p]$. LSPIA takes a weighted sum of all differences and computes the $(i,j)$th adjusting vector based on
\begin{align*}
{\bm \Delta}_{ij}^{(k)}= \widetilde{\mu} \sum_{h=0}^m \sum_{l=0}^p \mu_i(x_h) \mu_j(y_l){\bm R}_{hl}^{(k)}.
\end{align*}
Then, the new surface is generated by
\begin{align*}
  \mathrm{S}^{(k+1)}(x,y)=\sum_{i=0}^n \sum_{j=0}^n \mu_i(x) \mu_j(y) {\bm P}_{ij}^{(k+1)}
~ {\rm with} ~
{\bm P}_{ij}^{(k+1)} = {\bm P}_{ij}^{(k)} +  {\bm \Delta}_{ij}^{(k)}.
\end{align*}

It is a similar story to choose different ${\bm \Delta}_{ij}^{(k)}$ and obtain several specific LSPIA variants for surface fitting. The details here are omitted.
\vspace{.5ex}
\begin{remark}\label{remark:Landweber}
Let the collocation matrices of $\left\{\mu_i(x):i\in [n]\right\}$ on $\left\{x_h \right\}_{h=0}^m$  and $\left\{\mu_j(y):j\in [n]\right\}$ on $\left\{y_l \right\}_{l=0}^p$ respectively be
\begin{align*}
 {\bm A} = \begin{bmatrix} \mu_i(x_h)  \end{bmatrix}_{h=0,i=0}^{m,n} \in \Rc^{(m+1)\times (n+1)}
~and~
 {\bm B}^T = \begin{bmatrix}\mu_j(y_l)\end{bmatrix}_{l=0,j=0}^{p,n} \in \Rc^{(p+1)\times (n+1)}.
\end{align*}
We redemonstrate the forms of LSPIA iteration using linear algebra formulations as follows.
 \vspace{.5ex}
\begin{enumerate}[\quad\sffamily(1)]
       \setlength{\itemindent}{0cm}
       \item  The case of curves. Let the data points and control points be arranged respectively into
\begin{align*}
{\bm q}=\left[{\bm q}_{0}~{\bm q}_{1}~\cdots ~{\bm q}_{m}\right]^T
~and~
{\bm p}^{(k)}=\left[ {\bm p}_{0}^{(k)}~{\bm p}_{1}^{(k)}~\cdots~{\bm p}_{n}^{(k)}\right]^T
\end{align*}
for $k=0,1,2,\cdots$.  The LSPIA iterative process is expressed by
\begin{align}\label{LSPIA+ColumnIteration}
{\bm p}^{(k+1)} = {\bm p}^{(k)} + \widetilde{\mu} {\bm A}^T \left({\bm q} - {\bm A} {\bm p}^{(k)} \right).
\end{align}
       \item The case of surfaces. Let the data points and control points be arranged respectively into
\begin{align*}
\widetilde{{\bm Q}} = \left[{\bm Q}_{00}~\cdots~{\bm Q}_{m0}~{\bm Q}_{01}~\cdots~{\bm Q}_{m1}~ \cdots ~{\bm Q}_{0p}~\cdots~{\bm Q}_{mp}~\right]^T,
\end{align*}
and
\begin{align*}
\widetilde{{\bm P}}^{(k)} = \left[{\bm P}_{00}^{(k)}~\cdots~{\bm P}_{n0}^{(k)}~{\bm P}_{01}^{(k)}~\cdots~{\bm P}_{n1}^{(k)}~ \cdots ~{\bm P}_{0n}^{(k)}~\cdots~{\bm P}_{nn}^{(k)}~  \right]^T
\end{align*}
for $k=0,1,2,\cdots$.  The LSPIA iterative process corresponds to
\begin{align}
\widetilde{{\bm P}}^{(k+1)} = \widetilde{{\bm P}}^{(k)} + \widetilde{\mu} \left( {\bm B} \otimes {\bm A}^T \right) \left( \widetilde{{\bm Q}} - \left( {\bm B}^T \otimes {\bm A} \right)\widetilde{{\bm P}}^{(k)} \right).
\end{align}
\end{enumerate}
 \vspace{.5ex}
From algebraic aspects, LSPIA is equivalent to the Richardson method applied to a least-squares system \cite{Saad2000}.  In addition, formula  \eqref{LSPIA+ColumnIteration} is called Landweber's iteration in algebraic image reconstruction field, which is a row projection method \cite{51Landweber} and a special case of simultaneous iterative reconstruction technique \cite{ENH10}. It is also known as the gradient descent method in optimization and inverse problems; see, for example \cite{19JJ, 22RJ}.
\end{remark}

\vskip 1ex
Take the curve fitting situation as an example. Let ${\bm e}_i$ be the $i$th column of ${\bm I}_{n+1}$, where  ${\bm I}_{\ell}$ is the identity matrix with size $\ell$. By partitioning ${\bm A}$ into columns, i.e., ${\bm A}=[{\bm A}_{:,1}~{\bm A}_{:,2}~\cdots~{\bm A}_{:,n+1}]$, we know that
 \begin{align}\label{eq:LSPIA_ite}
{\bm p}^{(k+1)} = {\bm p}^{(k)} + \widetilde{\mu} \sum_{i=1}^{n+1} {\bm e}_{i} {\bm A}_{:,i}^T  \left( {\bm q} - {\bm A} {\bm p}^{(k)} \right).
\end{align}
Evidently, LSPIA needs to simultaneously adjust all ${\bm p}_{i}^{(k)}$ for $i\in[n]$. This discovery is heuristic. If we operate on part of ${\bm p}_{i}^{(k)}$, instead of all, we get the local LSPIA-type method.

\section{The RPIA method}\label{Sec:3}
In this section, we present the RPIA method  for curve and surface fittings. The process goes into detail on below.
\subsection{The case of curves}\label{sec:RPIA3.1}
We call $\left\{\mathrm{I}_i\right\}_{i=0}^{\imath}$ a partition of $[n]$ if $\mathrm{I}_i \cap \mathrm{I}_j = \emptyset$ for $i\neq j$ and $\cup_{i=0}^\imath \mathrm{I}_i = [n]$.  We first construct an initial curve
\begin{align*}
\mathrm{C}^{(0)}(x)=\sum_{i=0}^n \mu_i(x) {\bm p}_{i}^{(0)},~x\in [x_0,~x_m]
\end{align*}
and compute the $j$th difference according to
\begin{align*}
{\bm r}_{j}^{(0)}
= {\bm q}_{j} - \mathrm{C}^{(0)}(x_j)
\end{align*}
for $j\in [m]$. We randomly select $\imath_0\in [\imath]$ with probability
\begin{equation*}
  \bP\left({\rm Index} = {\imath_0} \right) =  \frac{\sum_{i_0 \in \mathrm{I}_{\imath_0}}  \sum_{j=0}^{m}\mu_{i_0}^2(x_j)}{\sum_{i=0}^{n}\sum_{j=0}^{m}\mu_i^2(x_j)}
\end{equation*}
and if $i_0 \in \mathrm{I}_{\imath_0} \subseteq [n]$, calculate the $i_0$th adjusting vector in keeping with
\begin{align*}
{\bm \delta}_{i_0}^{(0)} =
 \frac{\sum_{j=0}^m \mu_{i_0}(x_j){\bm r}_{j}^{(0)}}{\sum_{i_0 \in \mathrm{I}_{\imath_0}} \sum_{j=0}^{m}\mu_{i_0}^2(x_j)},
\end{align*}
otherwise ${\bm \delta}_{i_0}^{(0)}$ being zero. Then we update partial  control points in the light of
\begin{align*}
 {\bm p}_{\mathrm{I}_{\imath_0}}^{(1)} = {\bm p}_{\mathrm{I}_{\imath_0}}^{(0)} + {\bm \delta}_{\mathrm{I}_{\imath_0}}^{(0)}
\end{align*}
in which the rest of control points remain unchanged and ${\bm x}_{\mathrm{I}}$ denotes the subvector of ${\bm x}$ indexed by $\mathrm{I}$ for any vector ${\bm x}$.


Recursively, assuming that we have obtained the $k$th curve
\begin{align}\label{Curve+eq:C(x)_k}
\mathrm{C}^{(k)}(x)=\sum_{i=0}^n \mu_i(x) {\bm p}_{i}^{(k)},~x\in [x_0,~x_m]
\end{align}
for $k=0,1,2,\cdots$ and the $j$th difference
\begin{align}\label{Curve+eq:Update+rk}
{\bm r}_{j}^{(k)}
= {\bm q}_{j} - \mathrm{C}^{(k)}(x_j)
\end{align}
for $j\in [m]$. Let the index $\imath_k\in [\imath]$ be selected  with probability
\begin{equation}\label{Curve+eq:IndexChosenI}
  \bP\left({\rm Index} = {\imath_k} \right) =  \frac{\sum_{i_k \in \mathrm{I}_{\imath_k}}  \sum_{j=0}^{m}\mu_{i_k}^2(x_j)}{\sum_{i=0}^{n}\sum_{j=0}^{m}\mu_i^2(x_j)}.
\end{equation}
If $i_k \in \mathrm{I}_{\imath_k}\subseteq [n]$, the $i_k$th adjusting vector is computed by
\begin{align}\label{Curve+eq:Update+deltak}
{\bm \delta}_{i_k}^{(k)} =  \frac{\sum_{j=0}^m \mu_{i_k}(x_j) {\bm r}_{j}^{(k)}}{
 \sum_{i_k \in \mathrm{I}_{\imath_k}} \sum_{j=0}^{m}\mu_{i_k}^2(x_j) },
\end{align}
otherwise it is zero.
Then we update the control points in accordance with
\begin{equation}\label{Curve+eq:Update+pk}
 {\bm p}_{\mathrm{I}_{\imath_k}}^{(k+1)} = {\bm p}_{\mathrm{I}_{\imath_k}}^{(k)} + {\bm \delta}_{\mathrm{I}_{\imath_k}}^{(k)},
\end{equation}
while the other control points keep fixed. After that, the next curve is generated by
\begin{align*}
\mathrm{C}^{(k+1)}(x)=\sum_{i=0}^n  \mu_i(x) {\bm p}_{i}^{(k+1)}=\mathrm{C}^{(k)}(x) + \sum_{i_k \in \mathrm{I}_{\imath_k}}  \mu_{i_k}(x) {\bm \delta}_{i_k}^{(k)}.
\end{align*}
The RPIA for curve fitting is arranged in Algorithm \ref{Curve+alg:RPIA}.
\begin{algorithm}[htb]
\caption{ RPIA for curve fitting.}
\label{Curve+alg:RPIA}
\begin{algorithmic}[1]
\Require
Data points $\left\{{\bm  q}_{j} \right\}_{j=0}^{m}$,
the initial control points $\left\{{\bm p}_{i}^{(0)} \right\}_{i=0}^n$,
a real increasing sequence $\left\{x_j \right\}_{j=0}^m$,
an index set $\left\{\mathrm{I}_i\right\}_{i=0}^{\imath}$,
and the maximum iteration number $\ell$.
\Ensure
$\mathrm{C}^{(\ell)}(x)$.
\State {\bf for} $k=0,1,2,\cdots,\ell$ {\bf do}
\State \quad generate the blending curve $\mathrm{C}^{(k)}(x)$ as \eqref{Curve+eq:C(x)_k};
\State \quad calculate the differences as \eqref{Curve+eq:Update+rk};
\State \quad randomly pick the index $\imath_k$ as \eqref{Curve+eq:IndexChosenI};
\State \quad compute the adjusting vectors as \eqref{Curve+eq:Update+deltak};
\State \quad update the control points as \eqref{Curve+eq:Update+pk};
\State {\bf endfor}.
\end{algorithmic}
\end{algorithm}

\vskip .5ex
\begin{remark}\label{remark:VSlspia}
For a further insight into RPIA for curve fitting (Algorithm \ref{Curve+alg:RPIA}), we rewrite it into  matrix form. At the $k$th iteration,
\begin{align}\label{eq:RPIA_ite+curve}
 {\bm p}^{(k+1)} = {\bm p}^{(k)} + \widetilde{\mu}_k \sum_{i_k \in \mathrm{I}_{\imath_k}}{\bm e}_{i_k} {\bm A}_{:,i_k}^T \left({\bm q} - {\bm A} {\bm p}^{(k)}\right),
\end{align}
where $\widetilde{\mu}_k = 1/\|{\bm A}_{:,\mathrm{I}_{\imath_k}}\|_F^2$ and $\|\cdot \|_F$  represents the $F$-norm of a matrix. We note that in the sequel ${\bm M}_{:,\mathrm{I}}$ (resp. ${\bm M}^T_{\mathrm{J},:}$) denotes the column submatrix of ${\bm M}$ (resp., ${\bm M}^T$) indexed by $\mathrm{I}$ (resp., $\mathrm{J}$). Compared with formula \eqref{eq:LSPIA_ite}, it is clear that only ${\bm p}^{(k)}_{i_k}$ indexed by $\mathrm{I}_{\imath_k}$ is updated and the others are unchanged. The geometric meaning of RPIA is shown in Figure \ref{fig:RPIA-plot}.
 \begin{figure}[!h]
	\centering
	\includegraphics[width=4in,height=3in]{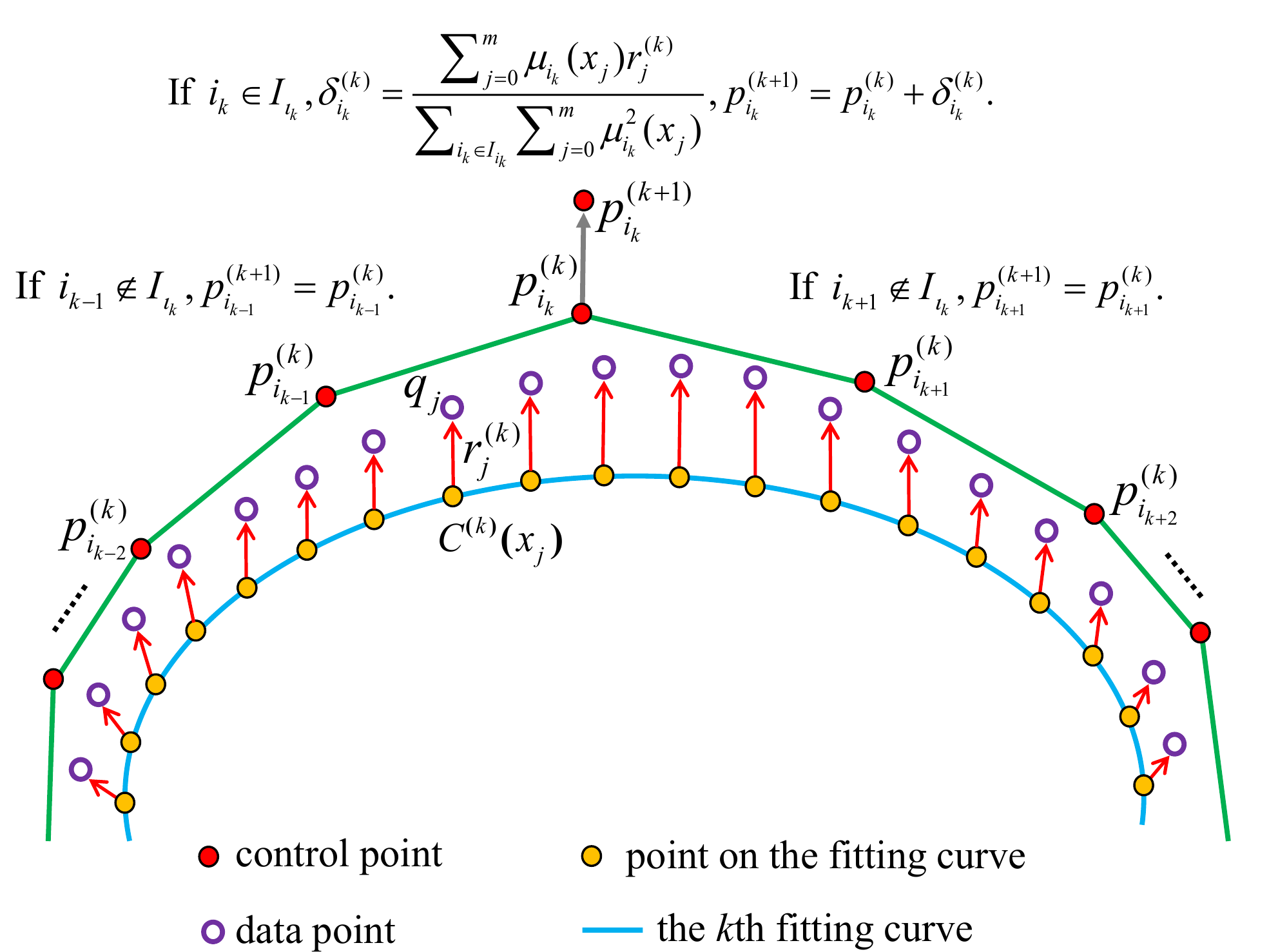}
	\caption{ Graphical illustration of RPIA for least-squares curve  fitting. }
	\label{fig:RPIA-plot}
\end{figure}
\end{remark}


\vskip .5ex
\begin{remark}\label{remark:RABCD}
 The algebraic interpretation of RPAI lies in the following aspect. Coordinate descent (CD, also known as column-oriented) method is popular in solving the general linear system  ${\bm M}{\bm x}={\bm g}$, where the coefficient matrix ${\bm M}\in\Rc^{(m+1) \times (n+1)}$, the right-hand side $\bg\in\Rc^{m+1}$, and the unknown vector $\bx\in\Rc^{n+1}$, due to its simplicity and numerical performance. The main idea of  CD is that, it operates on a column of ${\bm M}$ and chooses an unit coordinate direction as a search direction. The iterate is given  by
\begin{align*}
{\bm x}^{(k+1)}  = {\bm x}^{(k)} + \frac{1}{\bT{{\bm M}_{:,i}}^2} {\bm e}_{i} {\bm M}_{:,i}^T  \left({\bm g} - {\bm M}{\bm x}^{(k)} \right)
\end{align*}
for $k=0,1,2,\cdots$, where the index $i\in \{1,2, \cdots, n+1\}$ is chosen according to a well-defined criterion, such as a cyclic fashion or an appropriate probability distribution \cite{92LT}. To take advantage of parallel computation and further speed up the convergence, an extension of CD iterate is that, at each iteration multiple independent updates are computed and an average of them is used. Namely,
\begin{align*}
{\bm x}^{(k+1)}
& = {\bm x}^{(k)} + \sum_{i \in \mathrm{I}} \frac{p_i}{\bT{{\bm M}_{:,i}}^2} {\bm e}_{i} {\bm M}_{:,i}^T  \left({\bm g} - {\bm M}{\bm x}^{(k)} \right)\\
& = {\bm x}^{(k)} + \frac{1}{\bF{{\bm M}_{:,\mathrm{I}}}^2} \sum_{i \in \mathrm{I} } {\bm e}_{i} {\bm M}_{:,i}^T  \left({\bm g} - {\bm M}{\bm x}^{(k)} \right),
\end{align*}
where $p_i = \bT{{\bm M}_{:,i}}^2/\bF{{\bm M}_{:,\mathrm{I}}}^2$ with $\sum_{i \in \mathrm{I}}p_i=1 $ and the index set $\mathrm{I}$ is selected at random. It indicates that RPIA curve fitting equals to iteratively solve the linear system from Remark \ref{remark:VSlspia}. As far as we know, this randomized block CD method is new. For more discussions on the CD method, we refer to the works in \cite{10LL, 21NK, 15Wright, 20WX3} and the references therein.
\end{remark}

\subsection{The case of surfaces}\label{sec:RPIA3.2}
 We now turn to the case of surfaces. Let $\left\{\mathrm{J}_j\right\}_{j=0}^{\jmath}$ denote another partition of $[n]$.  We first construct an initial surface
\begin{align*}
\mathrm{S}^{(0)}(x,y)=\sum_{i=0}^n \sum_{j=0}^n \mu_i(x) \mu_j(y) {\bm P}_{ij}^{(0)},~x\in [x_0,~x_m],~y\in [y_0,~y_p]
\end{align*}
and compute the $(h,l)$th difference according to
\begin{align*}
{\bm R}_{hl}^{(0)} = {\bm Q}_{hl} - \mathrm{S}^{(0)}(x_h,y_l)
\end{align*}
for $h \in [m]$ and $l\in [p]$. Randomly selecting  $\imath_0\in [\imath]$  and  $\jmath_0\in [\jmath]$ with probabilities
\begin{align*}
  \bP\left({\rm Index} = {\imath_0} \right) =  \frac{\sum_{i_0 \in \mathrm{I}_{\imath_0}}  \sum_{h=0}^{m}\mu_{i_0}^2(x_h)}{\sum_{i=0}^{n}\sum_{h=0}^{m}\mu_i^2(x_h)}
~ {\rm and} ~
  \bP\left({\rm Index} = {\jmath_0} \right) =  \frac{\sum_{j_0 \in \mathrm{J}_{\jmath_0}}  \sum_{l=0}^{p}\mu_{j_0}^2(y_l)}{\sum_{j=0}^{n}\sum_{l=0}^{p}\mu_j^2(y_l)},
\end{align*}
respectively, we calculate the $(i_0, j_0)$th adjusting vector according to
\begin{align*}
{\bm \Delta}_{i_0, j_0}^{(0)} =  \frac{\sum_{h=0}^m \sum_{l=0}^p \mu_{i_0}(x_h)\mu_{j_0}(y_l){\bm R}_{hl}^{(0)}}{
\left( \sum_{i_0 \in \mathrm{I}_{\imath_0}} \sum_{h=0}^{m}\mu_{i_0}^2(x_h) \right)
\left( \sum_{j_0 \in \mathrm{J}_{\jmath_0}} \sum_{l=0}^{p}\mu_{j_0}^2(y_l) \right)}
\end{align*}
if $(i_0, j_0)\in (\mathrm{I}_{\imath_0}, \mathrm{J}_{\jmath_0})$ with $\mathrm{I}_{\imath_0}, \mathrm{J}_{\jmath_0}\subseteq [n]$, otherwise it is zero, and update the control points in the light of
\begin{align*}
 {\bm P}_{\mathrm{I}_{\imath_0}, \mathrm{J}_{\jmath_0}}^{(1)} = {\bm P}_{\mathrm{I}_{\imath_0}, \mathrm{J}_{\jmath_0}}^{(0)} + {\bm \Delta}_{\mathrm{I}_{\imath_0}, \mathrm{J}_{\jmath_0}}^{(0)},
\end{align*}
where the other control points remain unchanged.

Recursively, after  obtaining the $k$th surface
\begin{align}\label{Surface+eq:S(x,y)_k}
\mathrm{S}^{(k)}(x,y)=\sum_{i=0}^n \sum_{j=0}^n \mu_i(x) \mu_j(y) {\bm P}_{ij}^{(k+1)}
\end{align}
and the $(h,l)$th difference
\begin{align}\label{Surface+eq:Update+Rk}
{\bm R}_{hl}^{(k)}
= {\bm Q}_{hl} - \mathrm{S}^{(k)}(x_h,y_l)
\end{align}
for $h\in [m]$ and $l\in [p]$,  we randomly choose $\imath_k\in [\imath]$  and  $\jmath_k\in [\jmath]$ with probabilities
\begin{align}\label{Surface+eq:IndexChosenI}
  \bP\left({\rm Index} = {\imath_k} \right) =  \frac{\sum_{i_k \in \mathrm{I}_{\imath_k}}  \sum_{h=0}^{m}\mu_{i_k}^2(x_h)}{\sum_{i=0}^{n}\sum_{h=0}^{m}\mu_i^2(x_h)}
~ {\rm and} ~
  \bP\left({\rm Index} = {\jmath_k} \right) =  \frac{\sum_{j_k \in \mathrm{J}_{\jmath_k}}  \sum_{l=0}^{p}\mu_{j_k}^2(y_l)}{\sum_{j=0}^{n}\sum_{l=0}^{p}\mu_j^2(y_l)},
\end{align}
respectively, and compute the $(i_k, j_k)$th adjusting vector  based on
\begin{align}\label{Surface+eq:Update+Deltak}
{\bm \Delta}_{i_k, j_k}^{(k)} =  \frac{\sum_{h=0}^m \sum_{l=0}^p \mu_{i_k}(x_h)\mu_{j_k}(y_l){\bm R}_{hl}^{(k)}}{
\left( \sum_{i_k \in \mathrm{I}_{\imath_k}} \sum_{h=0}^{m}\mu_{i_k}^2(x_h) \right)
\left( \sum_{j_k \in \mathrm{J}_{\jmath_k}} \sum_{l=0}^{p}\mu_{j_k}^2(y_l) \right)}
\end{align}
if $(i_k, j_k) \in (\mathrm{I}_{\imath_k}, \mathrm{J}_{\jmath_k})$ with $\mathrm{I}_{\imath_0}, \mathrm{J}_{\jmath_0}\subseteq [n]$, otherwise it is zero. Then, the control points are updated by
\begin{equation}\label{Surface+eq:Update+Pk}
 {\bm P}_{\mathrm{I}_{\imath_k}, \mathrm{J}_{\jmath_k}}^{(k+1)} = {\bm P}_{\mathrm{I}_{\imath_k}, \mathrm{J}_{\jmath_k}}^{(k)} + {\bm \Delta}_{\mathrm{I}_{\imath_k}, \mathrm{J}_{\jmath_k}}^{(k)},
\end{equation}
where  the other control points retain fixed. With the above preparation, the next surface is generated by
\begin{align*}
\mathrm{S}^{(k+1)}(x,y)=\sum_{i=0}^n \sum_{j=0}^n \mu_i(x) \mu_j(y) {\bm P}_{ij}^{(k+1)}
=\mathrm{S}^{(k)}(x,y) +
\sum_{i_k \in \mathrm{I}_{\imath_k}}
\sum_{j_k \in \mathrm{J}_{\jmath_k}}
\mu_{i_k}(x) \mu_{j_k}(y){\bm \Delta}_{i_k, j_k}^{(k)}.
\end{align*}
The RPIA surface fitting is organized  in Algorithm \ref{Surface+alg:RPIA}.
\begin{algorithm}[htb]
\caption{ RPIA for surface fitting.}
\label{Surface+alg:RPIA}
\begin{algorithmic}[1]
\Require
Data points $\left\{{\bm Q}_{hl} \right\}_{h,l=0}^{m,p}$,
the initial control points $\left\{{\bm P}_{ij}^{(0)} \right\}_{i,j=0}^{n,n}$,
two real increasing sequences $\left\{x_h \right\}_{h=0}^m$ and $\left\{y_l \right\}_{l=0}^p$,
two  index sets $\left\{\mathrm{I}_i\right\}_{i=0}^{\imath}$ and $\left\{\mathrm{J}_j\right\}_{j=0}^{\jmath}$,
and the maximum iteration number $\ell$.
\Ensure
$\mathrm{S}^{(\ell)}(x,y)$.
\State {\bf for} $k=0,1,2,\cdots,\ell$ {\bf do}
\State \quad generate the blending surface $\mathrm{S}^{(k)}(x,y)$ as \eqref{Surface+eq:S(x,y)_k};
\State \quad calculate the differences as \eqref{Surface+eq:Update+Rk};
\State \quad randomly pick the indices $\imath_k$ and  $\jmath_k$ as \eqref{Surface+eq:IndexChosenI};
\State \quad compute the adjusting vectors as \eqref{Surface+eq:Update+Deltak};
\State \quad update  the control points as \eqref{Surface+eq:Update+Pk};
\State {\bf endfor}.
\end{algorithmic}
\end{algorithm}

\vskip .5ex
In Section \ref{sec:RPIA3.1}, we formulated a matrix expression for RPIA curve fitting. This approach carries over well to the case of surfaces. That is, we put the coordinates of each data point and  control point into a row partition  and directly apply Algorithm \ref{Curve+alg:RPIA} to a vectorized linear system. However, this poses a higher computational cost. We relieve this issue by taking an alternative tack. We write the iterative process in Algorithm \ref{Surface+alg:RPIA} in tensor form and equivalently turn to solve the matrix equations in the $x$-, $y$-, and $z$-axis directions, which is elaborated on below.

At $k$th iteration, assume that each data point and control point have the vector forms
\begin{align*}
 {\bm Q}_{hl}       = \left[ Q_{hl1}~Q_{hl2}~Q_{hl3}\right]^T\in \mathbb{R}^3
 ~{\rm and}~
 {\bm P}^{(k)}_{ij} = \left[ P^{(k)}_{ij1}~P^{(k)}_{ij2}~P^{(k)}_{ij3}\right]^T\in \mathbb{R}^3
\end{align*}
for $h\in [m]$, $l\in[p]$, $i,j\in [n]$, which respect the third-order tensor structures as follows.
\begin{align*}
 {\bm Q}= \left[Q_{hlt}^{} \right]_{h=0,l=0,t=1}^{m,p,3}~{\rm and}~
 {\bm P}^{(k)} = \left[P^{(k)}_{ijt} \right]_{i=0,j=0,t=1}^{n,n,3},
\end{align*}
with the $t$th frontal slice of ${\bm Q}$ and ${\bm P}^{(k)}$  being respectively defined by
\begin{align*}
{\bm Q}_{(t)} =
\begin{bmatrix}
Q_{00t}&Q_{01t}&\cdots&Q_{0pt}\\
Q_{10t}&Q_{11t}&\cdots&Q_{1pt}\\
\vdots&\vdots&\ddots&\vdots\\
Q_{m0t}&Q_{m1t}&\cdots&Q_{mpt}\\
\end{bmatrix}
~{\rm and}~
{\bm P}^{(k)}_{(t)} =
\begin{bmatrix}
P^{(k)}_{00t}&P^{(k)}_{01t}&\cdots&P^{(k)}_{0nt}\\
P^{(k)}_{10t}&P^{(k)}_{11t}&\cdots&P^{(k)}_{1nt}\\
\vdots&\vdots&\ddots&\vdots\\
P^{(k)}_{n0t}&P^{(k)}_{n1t}&\cdots&P^{(k)}_{nnt}\\
\end{bmatrix}
\end{align*}
for $t=1,2,3$. The recursion of \eqref{Surface+eq:Update+Pk} can be expressed as
\begin{align}\label{eq:RLSPIAupdate3}
{\bm P}^{(k+1)}_{(t)}    & = {\bm P}^{(k)}_{(t)} + {\bm \Delta}_{(t)}^{(k)},
\end{align}
where
\begin{align}\label{eq:RPIA_ite+surface}
{\bm{\Delta}}_{(t)}^{(k)}
= \frac{1}{\bF{{\bm A}_{:,\mathrm{I}_{\imath_k}}}^2  \bF{{\bm B}_{\mathrm{J}_{\jmath_k},:}}^2}
\sum_{i_k \in \mathrm{I}_{\imath_k}}\sum_{j_k \in \mathrm{J}_{\jmath_k}}
{\bm e}_{i_k} \left( {\bm A}_{:,i_k}^T \left( {\bm Q}_{(t)} - {\bm A}{\bm P}_{(t)}^{(k)}{\bm B} \right) {\bm B}_{j_k,:}^T \right) {\bm e}_{j_k}^T
\end{align}
and the indices  $ \imath_k \in[\imath]$ and $\jmath_k\in[\jmath]$ are respectively selected with probabilities
\begin{align*}
 \bP\left({\rm Index} = \imath_k \right) = \frac{\bF{{\bm A}_{:,\mathrm{I}_{\imath_k}}}^2}{\BF{{\bm A}}}
~{\rm and}~
 \bP\left({\rm Index} = \jmath_k \right) = \frac{\bF{{\bm B}_{\mathrm{J}_{\jmath_k},:} }^2}{\BF{{\bm B}}}.
\end{align*}
Formula \eqref{eq:RPIA_ite+surface} tells us that the  computation of ${\bm{\Delta}}_{(t)}^{(k)}$ is appropriate for parallel.

\vspace{.5ex}
\begin{remark}
Consider the matrix equation ${\bm M}{\bm X} {\bm N}={\bm G}$, where the coefficient matrices ${\bm M}\in\Rc^{(m+1) \times (n+1)}$ and ${\bm N}\in\Rc^{(n+1) \times (p+1)}$, the right-hand side ${\bm G}\in\Rc^{(m+1)\times (p+1)}$, and the unknown matrix ${\bm X} \in\Rc^{(n+1)\times (n+1)}$. At the $k$th iteration, the CD iterate is given  by
\begin{align*}
{\bm X}^{(k+1)}  = {\bm X}^{(k)} + \frac{1}{\bT{{\bm M}_{:,i}}^2\bT{{\bm N}_{j,:}}^2} {\bm e}_{i} \left( {\bm M}_{:,i}^T  \left({\bm G} - {\bm M}{\bm X}^{(k)}{\bm N} \right){\bm N}_{j,:}^T \right) {\bm e}_{j}^T.
\end{align*}
Similar to Remark \ref{remark:RABCD}, the randomized block version of the above iterate is
\begin{align*}
{\bm X}^{(k+1)}
= {\bm X}^{(k)} +
\sum_{i\in \mathrm{I} } \sum_{j\in \mathrm{J} }
\frac{p_i q_j}{\bT{{\bm M}_{:,i}}^2\bT{{\bm N}_{j,:}}^2} {\bm e}_{i} \left( {\bm M}_{:,i}^T  \left({\bm G} - {\bm M}{\bm X}^{(k)}{\bm N} \right){\bm N}_{j,:}^T \right) {\bm e}_{j}^T,
\end{align*}
where $p_i$ and $q_j$ are two probabilities with $\sum_{i \in \mathrm{I} } p_i = \sum_{j \in \mathrm{J} } q_j = 1$ and the index set $(\mathrm{I}, \mathrm{J})$ is selected at random. Formula \eqref{eq:RPIA_ite+surface} emerges if one takes
$p_i = \bT{{\bm M}_{:,i}}^2/\bF{{\bm M}_{:,\mathrm{I}}}^2$ and
$q_j = \bT{{\bm N}_{j,:}}^2/\bF{{\bm N}_{\mathrm{J}, :}}^2$.
This block iterative method is new to the best of our knowledge. Then, we say RPIA surface fitting is equivalent to solving three linear matrix equations from algebraic aspects.
\end{remark}

\vspace{.5ex}
\begin{remark}
We emphasize that, compared with LSPIA, RPIAs in Algorithms \ref{Curve+alg:RPIA} and \ref{Surface+alg:RPIA}  are more flexible and permit partial control points to be adjusted, not all of them. This is because only the control point ${\bm p}_{i}^{(k)}$ (resp., ${\bm P}_{ij}^{(k)}$) indexed by $\mathrm{I}_{\imath_k}$ (resp., $\mathrm{I}_{\imath_k}$ and $\mathrm{J}_{\jmath_k}$) is updated, and other control points remain unchanged, see formula \eqref{eq:RPIA_ite+curve} (resp., \eqref{eq:RPIA_ite+surface}). It implies that RPIAs are local.  This local format saves the computational resources significantly, especially when the number of data points is large. Later, this comparison will become much more apparent for larger test instances in the numerical section.
\end{remark}

\section{Convergence analyses of RPIA}\label{Sec:4}
In this section, we utilize matrix theory to analyze the convergence of RPIA for least-squares fitting.

\subsection{The case of curves}\label{subsec+caseOFcurve}
By the linear algebra theory, ${\bm p}^{\ast}$ is a least-squares solution of ${\bm A}{\bm p} = {\bm q}$ if and only if ${\bm A}^T{\bm A}{\bm p}^{\ast} = {\bm A}^T{\bm q}$.  Now we give the convergence analysis of Algorithm \ref{Curve+alg:RPIA} in the following theorem.

\vspace{.5ex}
\begin{theorem}\label{Curve+thm:RPIA}
Let $\left\{\mu_i(x):i \in[n] \right\}$ be  a blending basis sequence and ${\bm A}$ be the corresponding collocation matrix on the real increasing sequence $\left\{x_h \right\}_{h=0}^m$. Suppose that ${\bm A}$ has a full-column rank, when the number of data points is larger than that of control points, the fitting curve sequence, generated by RPIA (see Algorithm \ref{Curve+alg:RPIA}),   converges to the least-squares fitting solution in expectation.
\end{theorem}

\vspace{.5ex}
\begin{proof}
 Let us introduce an auxiliary intermediate ${\bm z}^{(k)} = {\bm A}\left({\bm p}^{(k)} - {\bm p}^{\ast}\right)$ for $k=0,1,2,\cdots$, where  ${\bm p}^{\ast}:=({\bm A}^T{\bm A})^{-1}{\bm A}^T{\bm q}$ is a least-squares solution.  Multiplying both sides by the transpose of the coefficient matrix yields
\begin{align*}
{\bm A}^T{\bm z}^{(k)} = {\bm A}^T{\bm A}{\bm p}^{(k)} - ({\bm A}^T{\bm A}) ({\bm A}^T{\bm A})^{-1}{\bm A}^T{\bm q} = -{\bm A}^T {\bm r}^{(k)}
\end{align*}
with ${\bm r}^{(k)}={\bm q} - {\bm A}{\bm p}^{(k)}$. Combined with formula \eqref{eq:RPIA_ite+curve}, it implies that
\begin{align*}
A{\bm \delta}^{(k)}
= A \left[ {\bm \delta}_{0}^{(k)}~{\bm \delta}_{1}^{(k)}~\cdots~{\bm \delta}_{n}^{(k)}\right]^T
= \frac{{\bm A}_{:,\mathrm{I}_{\imath_k}}{\bm A}_{:,\mathrm{I}_{\imath_k}}^T {\bm r}^{(k)}}{\bF{{\bm A}_{:,\mathrm{I}_{\imath_k}}}^2}
= -\frac{{\bm A}_{:,\mathrm{I}_{\imath_k}}{\bm A}_{:,\mathrm{I}_{\imath_k}}^T {\bm z}^{(k)}}{\bF{{\bm A}_{:,\mathrm{I}_{\imath_k}}}^2}
\end{align*}
and
\begin{align*}
{\bm z}^{(k+1)}
 =  {\bm A}\left({\bm p}^{(k)} + {\bm \delta}^{(k)} - {\bm p}^{\ast}\right)
 =  \left( {\bm I}_{m+1} -  \frac{{\bm A}_{:,\mathrm{I}_{\imath_k}}{\bm A}_{:,\mathrm{I}_{\imath_k}}^T }{\bF{{\bm A}_{:,\mathrm{I}_{\imath_k}}}^2} \right) {\bm z}^{(k)}.
 \end{align*}
Let $\bE_{k}^{\mathrm{I}}[\cdot] = \bE [\cdot|\imath_0,\imath_1,\cdots,\imath_{k-1}]$ denote the conditional expectation conditioned on the first $k$ iterations of Algorithm \ref{Curve+alg:RPIA},
 where $\imath_\ell$ means that the $\imath_\ell$th index set $\mathrm{I}_{\imath_\ell}$ is chosen.  By taking this conditional expectation, it yields that
\begin{align*}
 \bE_{k}^{\mathrm{I}}{ \left[  {\bm z}^{(k+1)}  \right] }
 & = {\bm z}^{(k)} - \bE_{k}^{\mathrm{I}}{ \left[ \frac{{\bm A}_{:,\mathrm{I}_{\imath_k}}{\bm A}_{:,\mathrm{I}_{\imath_k}}^T }{\bF{{\bm A}_{:,\mathrm{I}_{\imath_k}}}^2} \right] }  {\bm z}^{(k)}\\
 & = {\bm z}^{(k)} -   \sum_{i=0}^\imath
 \frac{\bF{{\bm A}_{:,\mathrm{I}_i}}^2}{\BF{{\bm A}}}
 \frac{{\bm A}_{:,\mathrm{I}_{\imath_k}}{\bm A}_{:,\mathrm{I}_{\imath_k}}^T }{\bF{{\bm A}_{:,\mathrm{I}_{\imath_k}}}^2}
 {\bm z}^{(k)} \\
 & = \left( {\bm I}_{m+1} -  \frac{{\bm A}{\bm A}^T}{\BF{{\bm A}}}  \right) {\bm z}^{(k)}.
\end{align*}
Based on the law of total expectation and unrolling the recurrence, it gives that
\begin{align*}
\bE{ \left[  {\bm z}^{(k+1)}  \right] }
= \left( {\bm I}_{m+1} -  \frac{{\bm A}{\bm A}^T}{\BF{{\bm A}}}  \right)  \bE{ \left[  {\bm z}^{(k)}  \right] }
=\cdots
= \left( {\bm I}_{m+1} -  \frac{{\bm A}{\bm A}^T}{\BF{{\bm A}}}  \right)^{k+1} {\bm z}^{(0)}.
\end{align*}
Multiplying left by ${\bm A}^T$, we have
\begin{align*}
{\bm A}^T \bE{ \left[  {\bm z}^{(k+1)}  \right] }
= \left( {\bm I}_{n+1} -  \frac{{\bm A}^T{\bm A}}{\BF{{\bm A}}}  \right)^{k+1} {\bm A}^T {\bm z}^{(0)}.
\end{align*}
Since ${\bm A}$ is full of column rank, ${\bm A}^T{\bm A}$ is positive
definite. It leads to
\begin{align*}
0 \leq \rho\left( {\bm I}_{n+1} - \frac{{\bm A}^T{\bm A}}{\BF{{\bm A}}} \right)<1,
\end{align*}
where $\rho\left({\bm M}\right)$  is the spectral radius of ${\bm M}$. Therefore,
\begin{align*}
  \lim_{k\rightarrow \infty} \left( {\bm I}_{n+1} -  \frac{{\bm A}^T{\bm A}}{\BF{{\bm A}}}  \right)^{k}
  = {\bm O}_{n+1},
\end{align*}
where ${\bm O}_{n+1}$ is the rank zero matrix with size $n+1$.  It follows that
\begin{align*}
 \bE{ \left[  {\bm p}^{(\infty)}  \right] }
  = \left( {\bm A}^T {\bm A} \right)^{-1} \left( {\bm A}^T \bE{ \left[  {\bm z}^{(\infty)} \right] } + {\bm A}^T{\bm A} {\bm p}^{\ast}   \right)
  = \left( {\bm A}^T {\bm A} \right)^{-1} {\bm A}^T{\bm A} {\bm p}^{\ast}
  = {\bm p}^{\ast}.
\end{align*}
From algebraic aspects, it is found that the RPIA limit curve is the least-squares fitting result to the  given data points.
\hfill
\end{proof}

\subsection{The case of surfaces}
Similar to the case of curves, the surface sequence, generated by Algorithm \ref{Surface+alg:RPIA}, is convergent in expectation to
the least-squares fitting result for the given data points. The result is stated as follows.

\vspace{.5ex}
\begin{theorem}\label{thm:randomized-LSPIA-surface}
Let $\left\{\mu_i(x):i \in[n] \right\}$ and $\left\{\mu_j(y):j\in[n]\right\}$ be two blending basis sequences, and ${\bm A}$ and ${\bm B}^T$ be the corresponding collocation matrices on the real increasing sequences $\left\{x_h \right\}_{h=0}^m$ and $\left\{y_l \right\}_{l=0}^p$, respectively. Suppose that ${\bm A}$ and ${\bm B}^T$ have a full-column rank, when the number of data points is larger than that of control points, the fitting surface sequence, generated by RPIA (see Algorithm \ref{Surface+alg:RPIA}), converges to the least-squares fitting solution in expectation.
\end{theorem}

\vspace{.5ex}
\begin{proof} As a preparatory step, let us introduce an auxiliary third-order tensor
\begin{align*}
 {\bm Z}^{(k)}= \left[Z^{(k)}_{ijt} \right]_{i=0,j=0,t=1}^{m,p,3},
\end{align*}
whose $t$th frontal slice is given by
\begin{align*}
{\bm Z}^{(k)}_{(t)} = {\bm A} \left({\bm P}_{(t)}^{(k)} - {\bm P}_{(t)}^{\ast}\right){\bm B} =
\begin{bmatrix}
Z_{00t}^{(k)}&Z_{01t}^{(k)}&\cdots&Z_{0pt}^{(k)}\\
Z_{10t}^{(k)}&Z_{11t}^{(k)}&\cdots&Z_{1pt}^{(k)}\\
\vdots&\vdots&\ddots&\vdots\\
Z_{m0t}^{(k)}&Z_{m1t}^{(k)}&\cdots&Z_{mpt}^{(k)}\\
\end{bmatrix}
\end{align*}
 with ${\bm P}_{(t)}^{\ast}=({\bm A}^T{\bm A})^{-1}{\bm A}^T {\bm Q}_{(t)}{\bm B}^T({\bm B}^T{\bm B})^{-1}$ for $t=1,2,3$ and $k=0,1,2,\cdots$. Left multiplication by ${\bm A}^T$ and right multiplication by ${\bm B}^T$ to ${\bm Z}_{(t)}^{(k)}$ lead to
\begin{equation*}
{\bm A}^T{\bm Z}_{(t)}^{(k)}{\bm B}^T = {\bm A}^T{\bm A}{\bm P}_{(t)}^{(k)}{\bm B}{\bm B}^T -  {\bm A}^T {\bm A}{\bm P}_{(t)}^{\ast} {\bm B}{\bm B}^T = -{\bm A}^T{\bm R}_{(t)}^{(k)} {\bm B}^T
\end{equation*}
with ${\bm R}_{(t)}^{(k)} = {\bm Q} - {\bm A}{\bm P}_{(t)}^{(k)}{\bm B}$. It follows that
\begin{align*}
{\bm Z}_{(t)}^{(k+1)}
& = {\bm A} \left({\bm P}_{(t)}^{(k+1)} - {\bm P}_{(t)}^{\ast}\right){\bm B}\\
& = {\bm A} \left({\bm P}_{(t)}^{(k)} - {\bm P}_{(t)}^{\ast} +  {\bm \Delta}_{(t)}^{(k)} \right){\bm B}\\
& = {\bm Z}^{(k)}_{(t)}  - \frac{{\bm A}_{:,\mathrm{I}_{\imath_k}} {\bm A}_{:,\mathrm{I}_{\imath_k}}^T {\bm Z}_{(t)}^{(k)} {\bm B}_{\mathrm{J}_{\jmath_k},:}^T {\bm B}_{\mathrm{J}_{\jmath_k},:}}{\bF{{\bm A}_{:,\mathrm{I}_{\imath_k}}}^2   \bF{{\bm B}_{\mathrm{J}_{\jmath_k},:} }^2},
\end{align*}
where the last equality is from
\begin{align*}
{\bm A}{\bm \Delta}_{(t)}^{(k)}{\bm B}
& = \frac{{\bm A}{\bm I}_{:,\mathrm{I}_{\imath_k}}{\bm A}_{:,\mathrm{I}_{\imath_k}}^T {\bm R}_{(t)}^{(k)} {\bm B}_{\mathrm{J}_{\jmath_k},:}^T {\bm I}_{\mathrm{J}_{\jmath_k},:}{\bm B}}{\bF{{\bm A}_{:,\mathrm{I}_{\imath_k}}}^2   \bF{{\bm B}_{\mathrm{J}_{\jmath_k},:} }^2}\\
& = \frac{{\bm A}_{:,\mathrm{I}_{\imath_k}}{\bm I}_{:,\mathrm{I}_{\imath_k}}^T {\bm A}^T {\bm R}_{(t)}^{(k)} {\bm B}^T {\bm I}_{\mathrm{J}_{\jmath_k},:}^T {\bm B}_{\mathrm{J}_{\jmath_k},:}}{\bF{{\bm A}_{:,\mathrm{I}_{\imath_k}}}^2   \bF{{\bm B}_{\mathrm{J}_{\jmath_k},:} }^2}\\
& = -\frac{{\bm A}_{:,\mathrm{I}_{\imath_k}} {\bm A}_{:,\mathrm{I}_{\imath_k}}^T {\bm Z}_{(t)}^{(k)} {\bm B}_{\mathrm{J}_{\jmath_k},:}^T {\bm B}_{\mathrm{J}_{\jmath_k},:}}{\bF{{\bm A}_{:,\mathrm{I}_{\imath_k}}}^2   \bF{{\bm B}_{\mathrm{J}_{\jmath_k},:} }^2}.
\end{align*}
Let
$\bE_{k}^{\mathrm{I},\mathrm{J}}[\cdot] = \bE [\cdot|\jmath_0,\imath_0,\cdots,\jmath_{k-1},\imath_{k-1}]$
 denote the conditional expectation conditioned on the first $k$ iterations of Algorithm \ref{Surface+alg:RPIA},
 where $\imath_\ell$ and  $\jmath_{\ell}$ mean that the $\imath_\ell$th index set $\mathrm{I}_{\imath_\ell}$ and  the $\jmath_{\ell}$th index set $ \mathrm{J}_{\jmath_\ell}$ are chosen as
$\bE_{k}^{\mathrm{J}}\left[\cdot\right]= \bE \left[\cdot|\jmath_0,\imath_0,\cdots,\jmath_{k-1},\imath_{k-1},\imath_k\right]$
 and
$\bE_{k}^{\mathrm{I}}\left[\cdot\right]= \bE \left[\cdot|\jmath_0,\imath_0,\cdots,\jmath_{k-1},\imath_{k-1},\jmath_k\right]$,
respectively.  We have
\begin{align*}
 \bE_{k}^{\mathrm{I}, \mathrm{J}}{ \left[ {\bm Z}_{(t)}^{(k+1)} \right] }
 & = {\bm Z}_{(t)}^{(k)} -  \bE_{k}^{\mathrm{I}, \mathrm{J}}{ \left[ \frac{{\bm A}_{:,\mathrm{I}_{\imath_k}} {\bm A}_{:,\mathrm{I}_{\imath_k}}^T {\bm Z}_{(t)}^{(k)} {\bm B}_{\mathrm{J}_{\jmath_k},:}^T {\bm B}_{\mathrm{J}_{\jmath_k},:}}{\bF{{\bm A}_{:,\mathrm{I}_{\imath_k}}}^2   \bF{{\bm B}_{\mathrm{J}_{\jmath_k},:} }^2} \right] }\\
 & = {\bm Z}^{(k)}_{(t)} -
 \sum_{i=0}^\imath \sum_{j=0}^\jmath
 \frac{\bF{{\bm A}_{:,\mathrm{I}_i}}^2}{\BF{{\bm A}}}
 \frac{\bF{{\bm B}_{\mathrm{J}_j,:} }^2}{\BF{{\bm B}}}
 \frac{{\bm A}_{:,\mathrm{I}_{\imath_k}} {\bm A}_{:,\mathrm{I}_{\imath_k}}^T {\bm Z}^{(k)}_{(t)} {\bm B}_{\mathrm{J}_{\jmath_k},:}^T {\bm B}_{\mathrm{J}_{\jmath_k},:}}{\bF{{\bm A}_{:,\mathrm{I}_{\imath_k}}}^2   \bF{{\bm B}_{\mathrm{J}_{\jmath_k},:} }^2}\\
 & = {\bm Z}^{(k)}_{(t)} -  \frac{{\bm A}{\bm A}^T {\bm Z}_{(t)}^{(k)} {\bm B}^T{\bm B} }{\BF{{\bm A}}\BF{{\bm B}}}.
\end{align*}
By the law of total expectation, it yields that
\begin{align*}
 \bE { \left[ {\bm Z}_{(t)}^{(k+1)} \right] } =
 \bE { \left[ {\bm Z}_{(t)}^{(k)  } \right] }  -  \frac{{\bm A}{\bm A}^T \bE { \left[ {\bm Z}_{(t)}^{(k)  } \right] } {\bm B}^T{\bm B} }{\BF{{\bm A}}\BF{{\bm B}}}.
\end{align*}
Let ${\bm Z}^{(k)}_{(t)}$ be arranged into a row partition, i.e.,
\begin{align*}
\widetilde{\bm z}_{(t)}^{(k)}
 = \left[Z_{00t}^{(k)}~\cdots~Z_{m0t}^{(k)}~Z_{01t}^{(k)}~\cdots~Z_{m1t}^{(k)}~ \cdots ~Z_{0pt}^{(k)}~\cdots~Z_{mpt}^{(k)}~  \right]^T
\end{align*}
 for $t=1,2,3$ and $k=0,1,2\cdots$. It achieves that
\begin{align*}
 \bE { \left[ \widetilde{\bm z}_{(t)}^{(k+1)} \right] }
 & =
 \left({\bm I}_{(m+1)(p+1)} -  \left( \frac{{\bm B}^T{\bm B}}{\BF{{\bm B}}}\right) \otimes\left(\frac{{\bm A} {\bm A}^T}{\BF{{\bm A}}}\right)   \right)\bE { \left[ \widetilde{\bm z}_t^{(k)} \right] }\\
 & \qquad \qquad \qquad \qquad \qquad \vdots\\
 & = \left({\bm I}_{(m+1)(p+1)} -  \left( \frac{{\bm B}^T{\bm B}}{\BF{{\bm B}}}\right) \otimes\left(\frac{{\bm A} {\bm A}^T}{\BF{{\bm A}}}\right)\right)^{k+1}\widetilde{\bm z}_{(t)}^{(0)}.
\end{align*}
Multiplying left by ${\bm B}\otimes {\bm A}^T$, we have
\begin{align*}
\left( {\bm B}\otimes {\bm A}^T \right)\bE { \left[ \widetilde{\bm z}_{(t)}^{(k+1)} \right] }
= \left({\bm I}_{(n+1)^2} -  \left( \frac{{\bm B}{\bm B}^T}{\BF{{\bm B}}}\right) \otimes\left(\frac{{\bm A}^T {\bm A}}{\BF{{\bm A}}}\right)\right)^{k+1} \left( {\bm B}\otimes {\bm A}^T \right)\widetilde{\bm z}_{(t)}^{(0)}.
\end{align*}
Since ${\bm A}$ and ${\bm B}^T$ are full of column rank, ${\bm A}^T{\bm A}$ and ${\bm B}{\bm B}^T$ are positive
definite. It leads to
\begin{align*}
0 \leq \rho\left({\bm I}_{(n+1)^2} -  \left( \frac{{\bm B}{\bm B}^T}{\BF{{\bm B}}}\right) \otimes\left(\frac{{\bm A}^T {\bm A}}{\BF{{\bm A}}}\right)\right)<1.
\end{align*}
Therefore,
\begin{align*}
  \lim_{k\rightarrow \infty} \left({\bm I}_{(n+1)^2} -  \left( \frac{{\bm B}{\bm B}^T}{\BF{{\bm B}}}\right) \otimes\left(\frac{{\bm A}^T {\bm A}}{\BF{{\bm A}}}\right)\right)^{k}
  = {\bm O}_{(n+1)^2}.
\end{align*}
It follows that
\begin{align*}
 \bE{ \left[  {\bm  P}_{(t)}^{(\infty)}  \right] }
  = \left( {\bm A}^T {\bm A} \right)^{-1} \left( {\bm A}^T \bE{ \left[  {\bm Z}_{(t)}^{(\infty)} \right] }{\bm B}^T  + {\bm A}^T{\bm A} {\bm P}_{(t)}^{\ast}{\bm B} {\bm B}^T   \right)
  \left( {\bm B} {\bm B}^T \right)^{-1}
  = {\bm P}_{(t)}^{\ast},
\end{align*}
which indicates that the sequence of surfaces, generated by RPIA, converges to the least-squares fitting result in expectation.
\hfill
\end{proof}

\section{Numerical experiments}\label{Sec:5}
 In this section,  we give several representative examples and perform the RPIA method for curve and surface fittings. The cubic B-spline basis is used because of its simplicity and wide range of applications in computer-aided design, see \cite{14DL,19EL,20HW}.

 In Algorithms \ref{Curve+alg:RPIA} and \ref{Surface+alg:RPIA}, suppose that the subsets $\left\{ \mathrm{I}_{i} \right\}_{i=0}^{\imath}$ and $\left\{ \mathrm{J}_{j} \right\}_{j=0}^{\jmath}$ have the same size $\tau$, i.e., $|\mathrm{I}_{i}|=|\mathrm{J}_{j}|=\tau$ for any $i\in [\imath]$ and $j\in [\jmath]$. To be specific, we consider the two following partitions.
 \begin{align*}
   & \mathrm{I}_{i} = \mathrm{J}_{j}=
   \left\{i\tau+1, i\tau+2,\cdots, i \tau + \tau \right\},
   && i=j = 0,1,2,\cdots,\imath-1,\\
   & \mathrm{I}_{\imath} = \mathrm{J}_{\jmath} =
   \left\{\imath\tau+1, \imath\tau+2,\cdots, n \right\},
   &&|\mathrm{I}_{\imath}| = |\mathrm{J}_{\jmath}| \leq \tau.
 \end{align*}
The symbol RPIA($\tau$) represents the RPIA method having the block size $\tau$. The discrete sampling is realized by applying MATLAB built-in function, e.g., {\sf  randsample}. We repeatedly run RPIA $30$ times and take the arithmetic mean of the results.  To make the implementation of RPIA more efficient, we try to avoid using for-loop structure as far as possible at each iteration.

We compare the performance of our method with LSPIA \cite{14DL}, SLSPIA \cite{19EL}, and MLSPIA \cite{20HW} in terms of iteration number (denoted as IT), computing time in seconds (denoted as CPU), and relative fitting error respectively defined by
\begin{equation*}
 E_k = \frac{\sum_{i=0}^n \BT{ \sum_{j=0}^m \mu_i(x_j){\bm r}^{(k)}_{j}}}{\sum_{i=0}^n \BT{ \sum_{j=0}^m \mu_i(x_j){\bm r}^{(0)}_{j}}}
~{\rm and}~
 E_k = \frac{\sum_{i=0}^n \sum_{j=0}^n \BT{ \sum_{h=0}^m \sum_{l=0}^p\mu_i(x_h)\mu_j(y_l){\bm R}^{(k)}_{hl}}}
 {\sum_{i=0}^n \sum_{j=0}^n \BT{ \sum_{h=0}^m \sum_{l=0}^p\mu_i(x_h)\mu_j(y_l){\bm R}^{(0)}_{hl}}}
\end{equation*}
for curve and surface cases when $k=0,1,2,\cdots$. The experiments are terminated once $ E_{k}$
is less than $10^{-6}$ or IT exceeds $10^{4}$ and let $E_k$ be $E_\infty$ \cite{14DL}.

As shown in  \cite[Section 3.2]{14DL}, \cite[Remark 2]{19EL}, and \cite[Theorems 6--7]{20HW}, the practical methods for selecting the appropriate weights appeared in LSPIA, SLSPIA, and MLSPIA are given by
\begin{equation*}
\widetilde{\mu}  = \frac{2}{\max_{i\in[n]}\left\{ \widetilde{c}_i\right\}},
~
\widehat{\mu} = \frac{2}{\max_{i\in[n]}\left\{ \widehat{c}_i  \right\}},
~
\omega =\gamma = \frac{4\sigma_1 \sigma_r}{(\sigma_1 + \sigma_r)^2},
~{\rm and}~
v = \frac{1}{\sigma_1\sigma_r},
\end{equation*}
where $\widetilde{c}_i$ and $\widehat{c}_i$ are the sums of the $i$th row elements of matrices ${\bm Y}^T{\bm Y}$ and $({\bm Y}^T{\bm Y})^2$, respectively; $\sigma_1$ and $\sigma_r$ are the largest and the smallest singular values of ${\bm Y}$, respectively; ${\bm Y}={\bm A}$ in the case of curves and ${\bm Y}={\bm B}^T \otimes {\bm A}$ in the case of surfaces. Note that all the singular values are computed via MATLAB function, e.g., {\sf svd}. We execute MLSPIA without explicitly forming $\sigma_1$ and $\sigma_r$.

The data points are from \cite{14DL}, \cite{20HW}, and the collection of various topics in geometry (available from \verb"http://paulbourke.net/geometry/.").

\subsection{Curve fitting}\label{subsec:NE1}
\begin{figure}[!htb]
\centering
	\subfigure[Example \ref{ex:RPIA1}]{
	\includegraphics[width=0.4\textwidth]{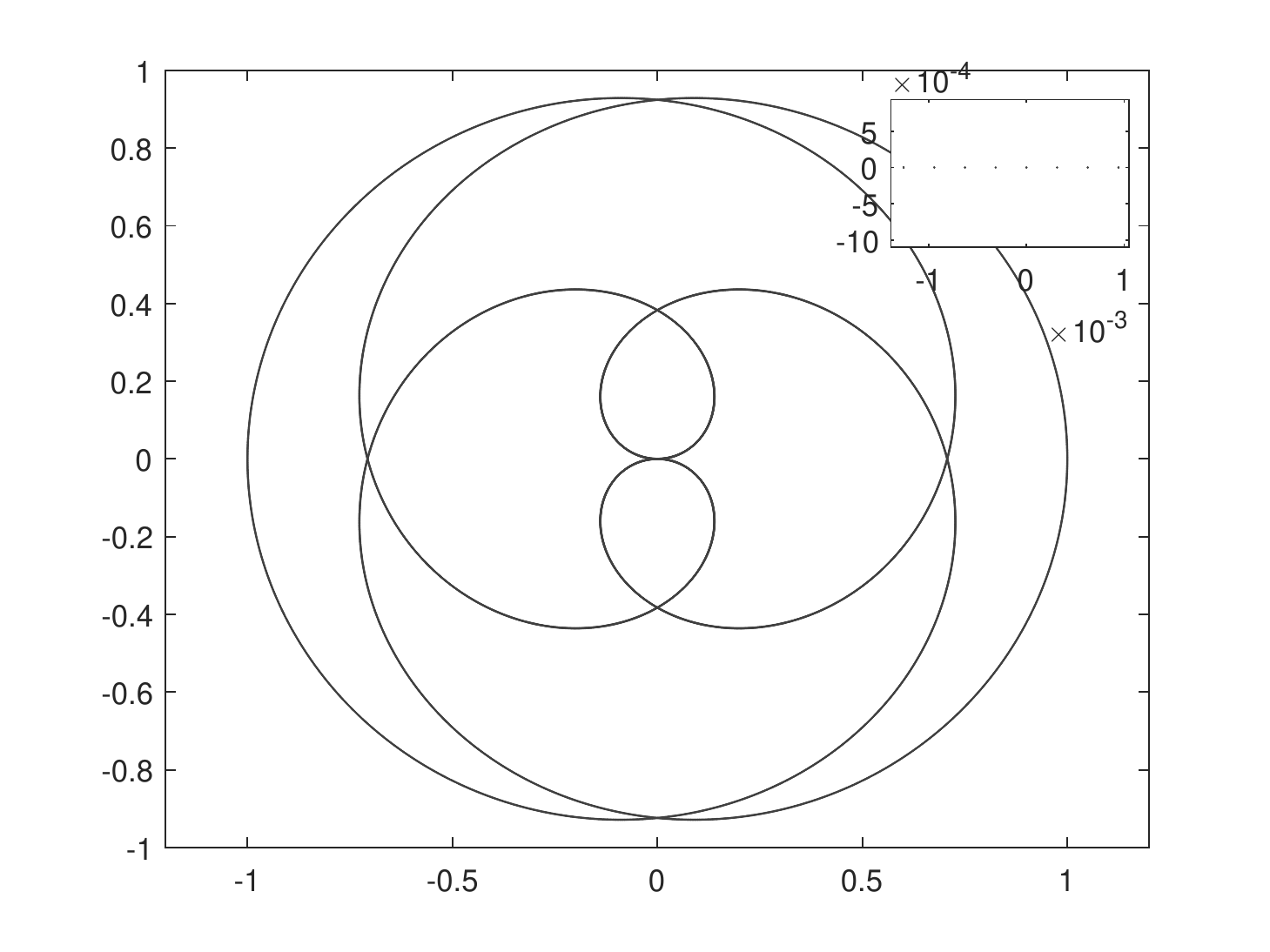}}
	\subfigure[Example \ref{ex:RPIA2}]{
	\includegraphics[width=0.4\textwidth]{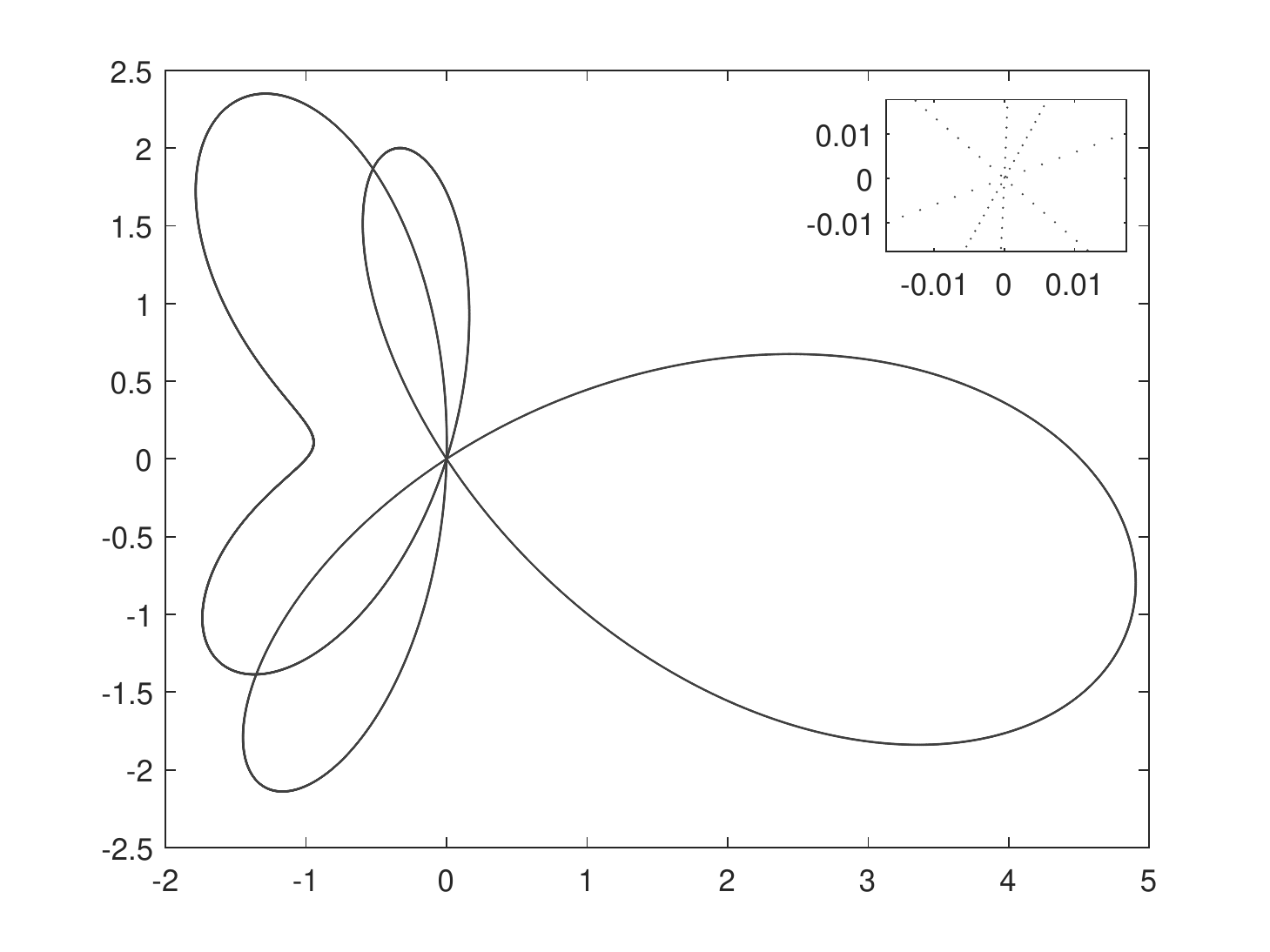}}
 	\subfigure[Example \ref{ex:RPIA3}]{
	\includegraphics[width=0.4\textwidth]{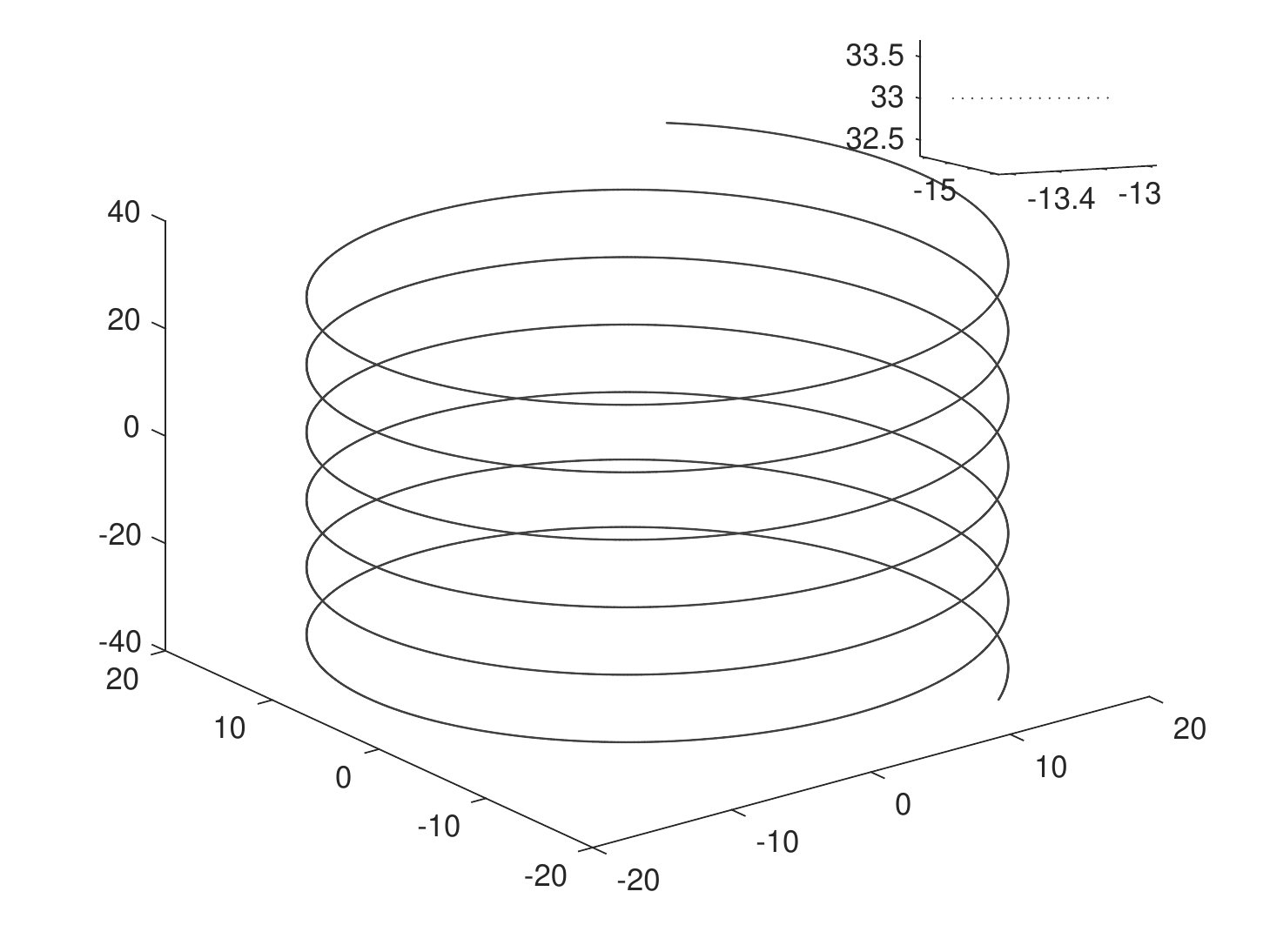}}
 	\subfigure[Example \ref{ex:RPIA4}]{
	\includegraphics[width=0.4\textwidth]{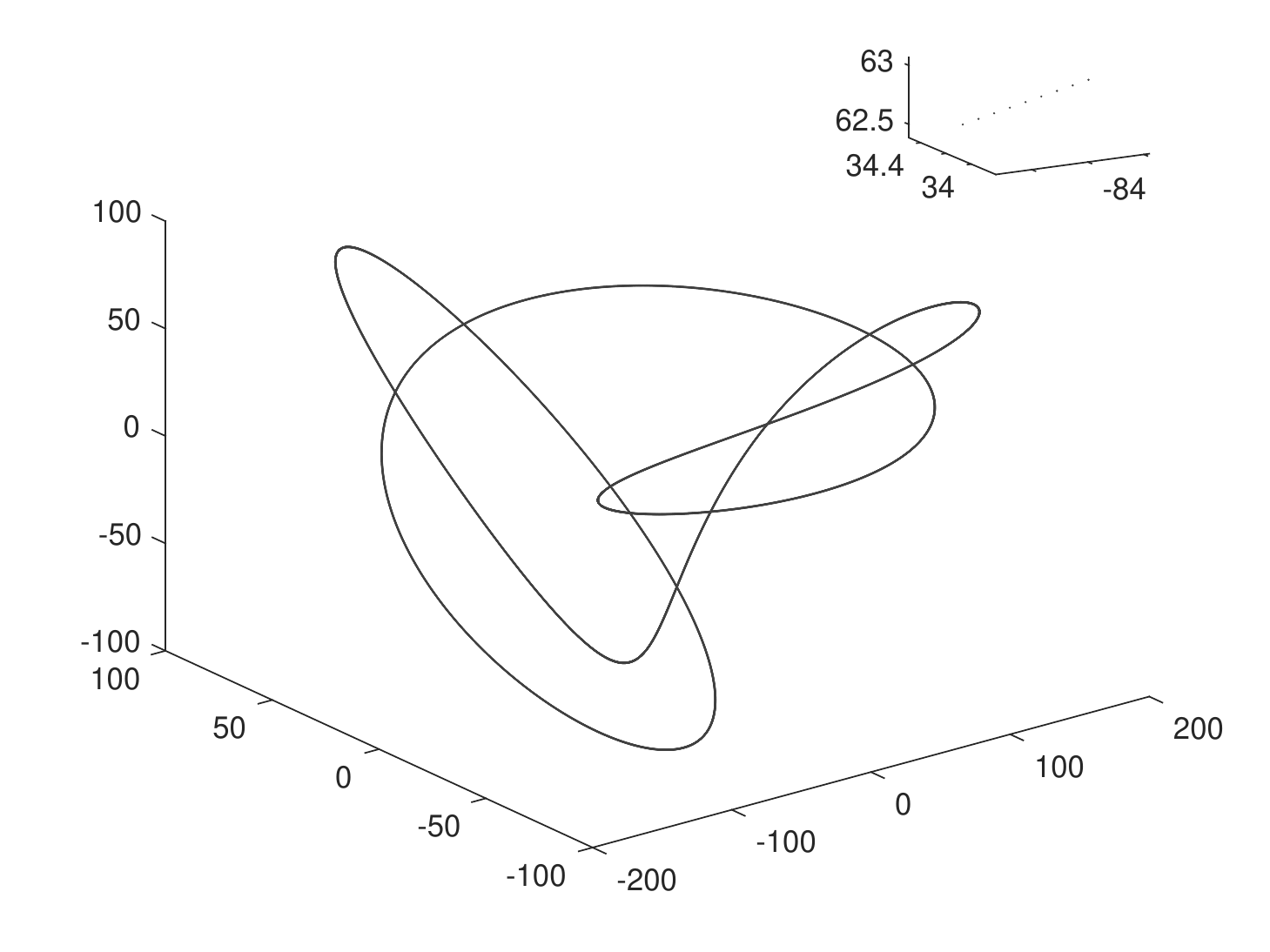}}
\caption{The data point sets to be fitted in Examples \ref{ex:RPIA1} (a), \ref{ex:RPIA2} (b), \ref{ex:RPIA2} (c), and \ref{ex:RPIA4} (d) with $m=20000$.}
\label{fig:curve+initial}
\end{figure}

The implementation details of RPIA for curve fitting are arranged as follows. We assign the parameters  $\left\{x_j \right\}_{j=0}^m$ for $\left\{{\bm q}_j\right\}_{j=0}^m$ according to the normalized accumulated chord  parameterization  method, i.e.,
\begin{equation*}
 x_0 = 0,~x_j = x_{j-1}+\frac{\bT{{\bm q}_j-{\bm q}_{j-1}}}{\sum_{s=1}^m \bT{{\bm q}_s-{\bm q}_{s-1}}},
~{\rm and}~
 x_m =1
\end{equation*}
for $j=1,2,\cdots,m$. The knot vector of cubic B-spline basis is defined by
\begin{equation*}
  \left\{
  0,~0,~0,~0,
  ~\bar{x}_4,~\bar{x}_5,\cdots,\bar{x}_n,
  ~1,~1,~1,~1
  \right\},
\end{equation*}
  with $\bar{x}_{j+3}=(1-\alpha)x_{i-1} + \alpha x_i$, $i = \lfloor jd \rfloor,~\alpha=jd-i$, $d = (m+1)/(n-2)$ for $j=1,2,\cdots,n-3$, and the notation $\lfloor \cdot \rfloor$ is the greatest integer function. For $i=1,2,\cdots,n-1$, the initial control points are selected by
 ${\bm p}_i^{(0)} = {\bm q}_{f_1(i)}$ with $f_1(0)= 0$, $f_1(n)=m$, and $f_1(i)=\lfloor  mi/n  \rfloor$ for $i=1,2,\cdots,n-1$, which is also described in the equation (23) of Deng and Lin \cite{14DL}.  Four point sets are considered and shown in Figure \ref{fig:curve+initial}.

\begin{example}\label{ex:RPIA1}
  $m+1$ data points sampled uniformly from a rose-type curve, whose polar coordinate equation is
  \begin{align*}
    {\bm r} =\sin(\theta/4) \quad (0 \leq \theta \leq 8\pi).
  \end{align*}
\end{example}

\begin{example}\label{ex:RPIA2}
  $m+1$ data points sampled uniformly from a blob-shaped curve, whose polar coordinate equation is
  \begin{align*}
    {\bm r} = 1+2\cos(2\theta+1/2) + 2\cos(3\theta+1/2) \quad (0 \leq \theta \leq 2\pi).
  \end{align*}
\end{example}

\begin{example}\label{ex:RPIA3}
  $m+1$ data points sampled from a helix curve, whose coordinates are given by
\begin{equation*}
\left\{ \begin{array}{l}
\bx = 10\cos( t\pi/3), \vspace{1ex}\\
\by = 10\sin(t\pi/3), \vspace{1ex}\\
\bz = t\pi/3~(-10\pi \leq t \leq 10\pi).
\end{array}\right .
\end{equation*}
\end{example}

\begin{example}\label{ex:RPIA4}
  $m+1$ data points sampled from a granny knot curve, whose coordinates are given by
\begin{equation*}
\left\{ \begin{array}{l}
\bx = -22\cos(t) - 128 \sin(t) - 44 \cos(3t) - 78 \sin(3t), \vspace{1ex}\\
\by = -10 \cos(2t) - 27 \sin(2t) + 38 \cos(4t) + 46 \sin(4t), \vspace{1ex}\\
\bz = 70 \cos(3t) - 40 \sin(3t)~ (0 \leq t \leq 2\pi).
\end{array}\right .
\end{equation*}
\end{example}

\begin{figure}[!htb]
\centering
 	\subfigure[Curve by RPIA(5)]{
		\includegraphics[width=0.4\textwidth]{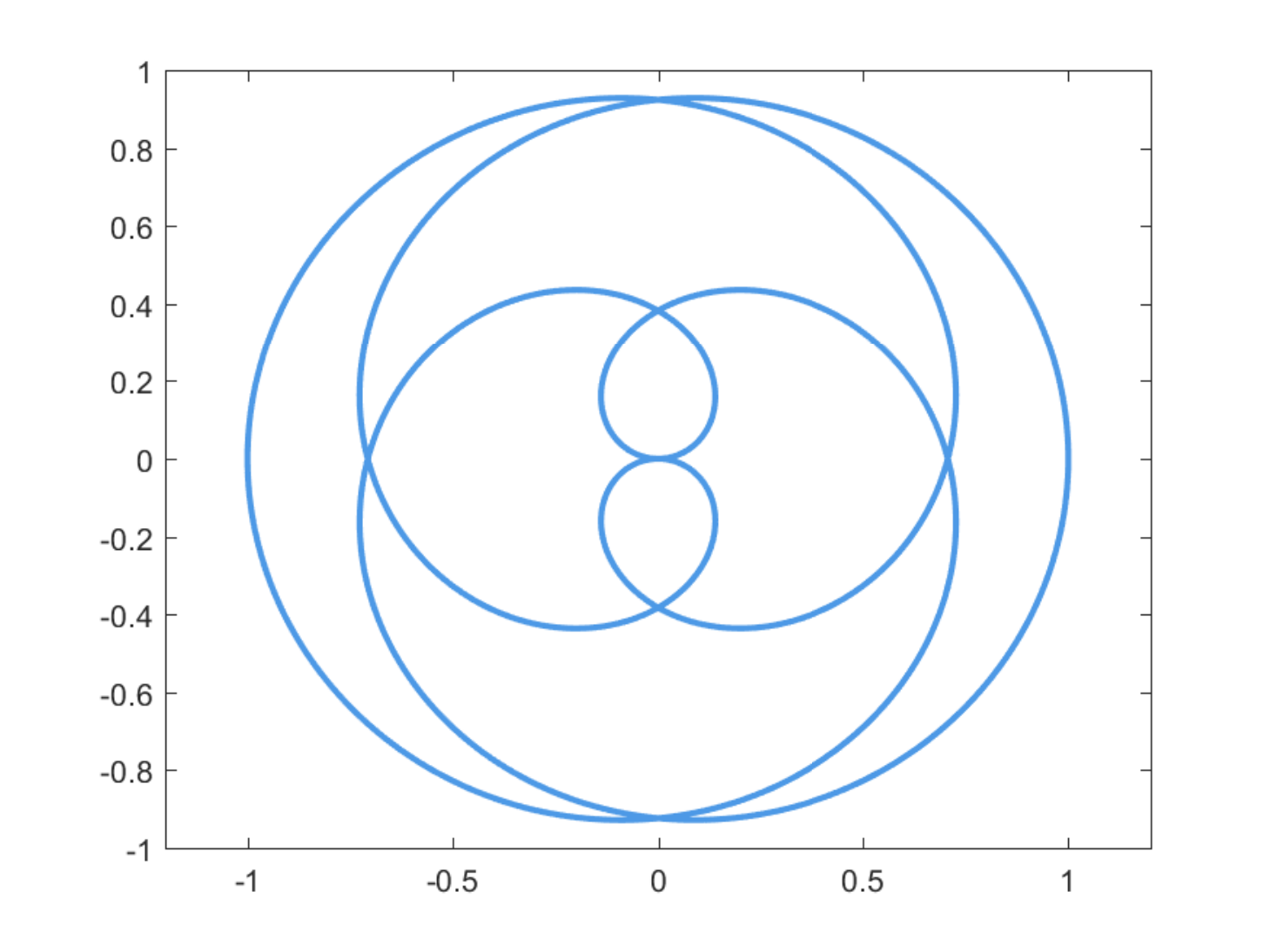}}
	\subfigure[Curve by RPIA(10)]{
	    \includegraphics[width=0.4\textwidth]{Ex1_Curve_m20000_n500-eps-converted-to.pdf}}
\caption{The cubic B-spline fitting curves given by RPIA with $m=20000$ and $n= 500$ for Example \ref{ex:RPIA1}.}
\label{fig:ex1-RPIA1+Result}
\end{figure}

\begin{figure}[!htb]
\centering
 	\subfigure[Curve by RPIA(5)]{
		\includegraphics[width=0.4\textwidth]{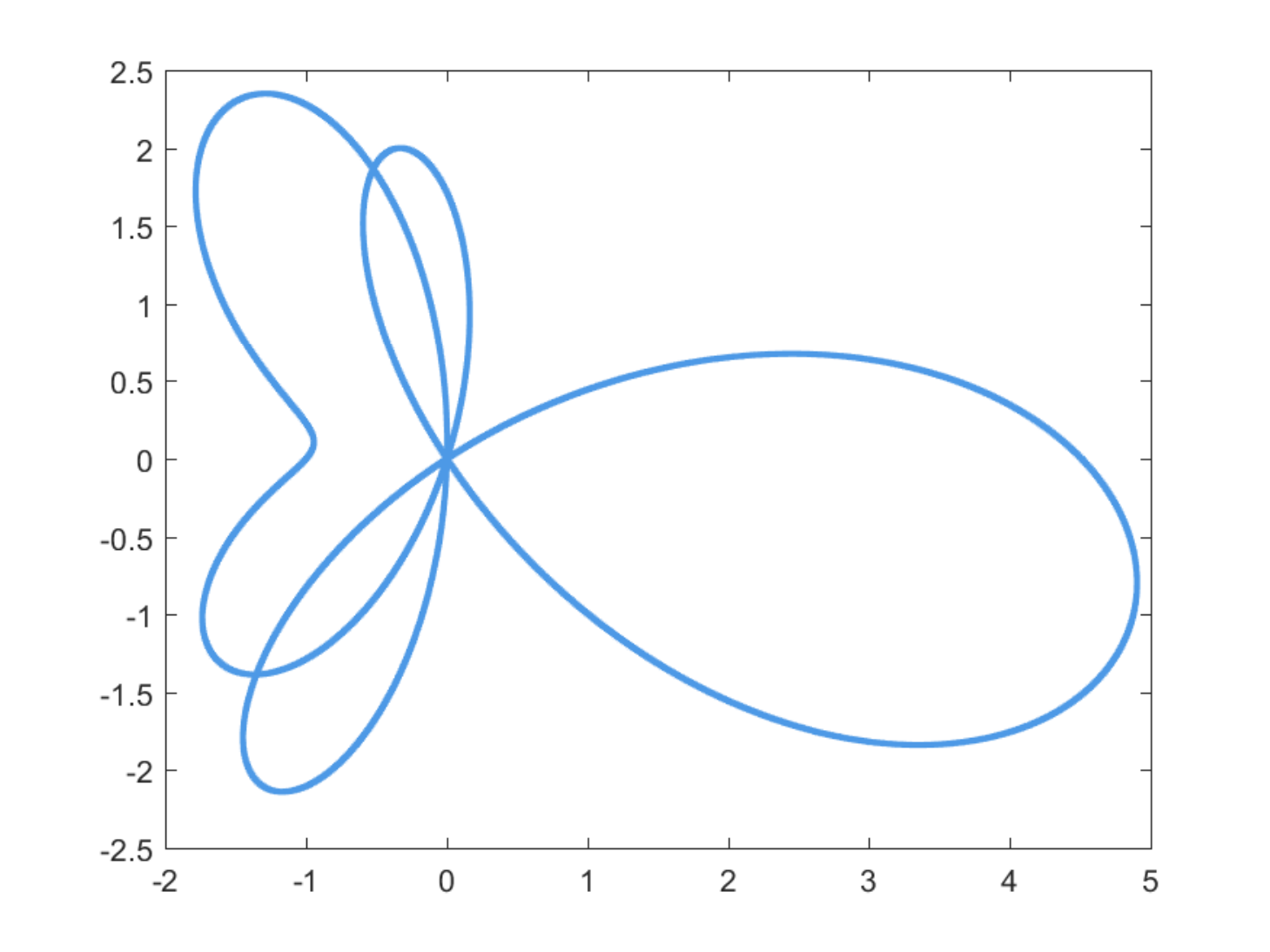}}
	\subfigure[Curve by RPIA(10)]{
	    \includegraphics[width=0.4\textwidth]{Ex2_Curve_m20000_n500-eps-converted-to.pdf}}
\caption{The cubic B-spline fitting curves given by RPIA with $m=20000$ and $n= 500$ for Example \ref{ex:RPIA2}.}
\label{fig:ex2-RPIA1+Result}
\end{figure}

\begin{figure}[!htb]
\centering
 	\subfigure[Curve by RPIA(5)]{
		\includegraphics[width=0.4\textwidth]{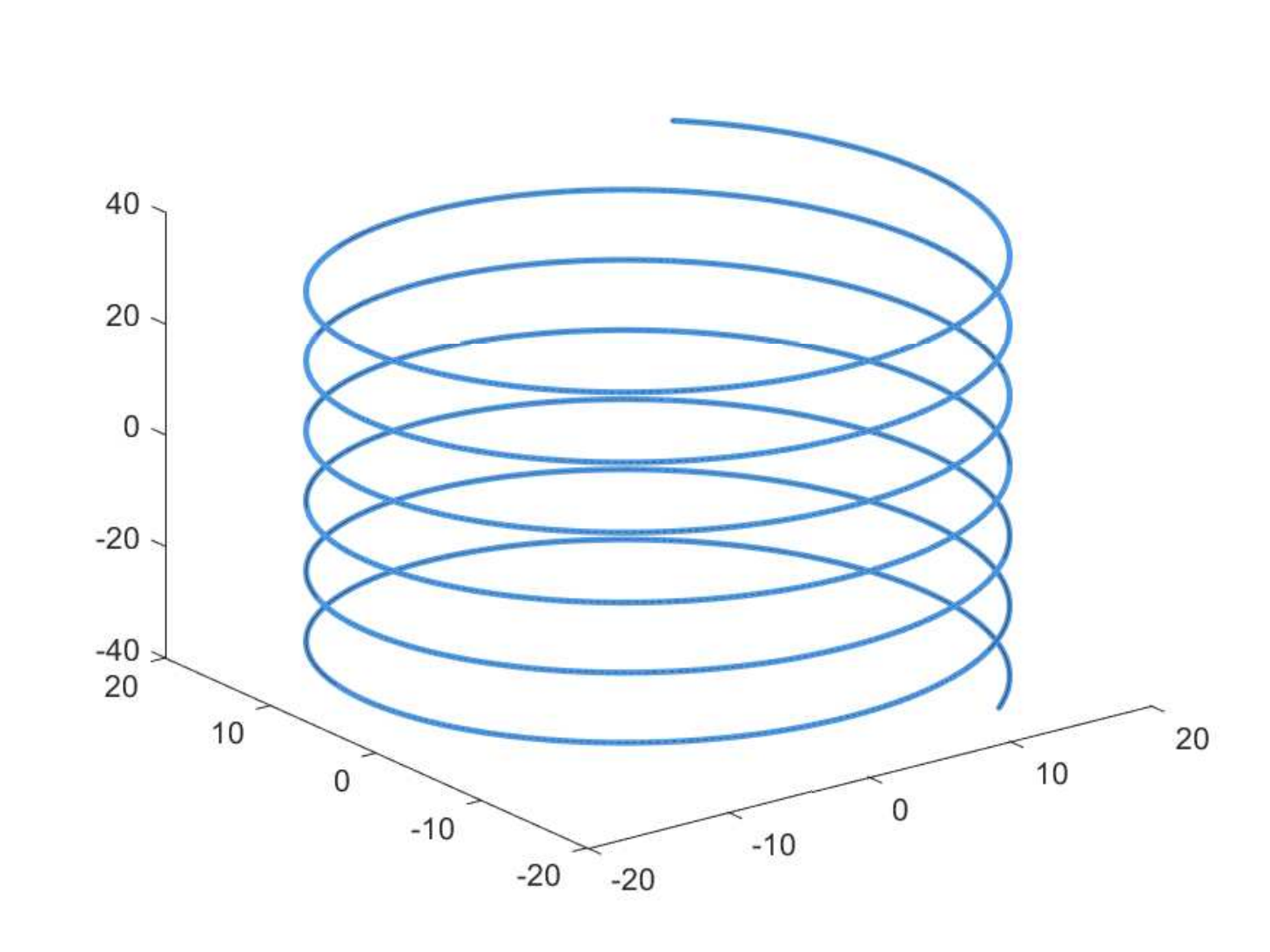}}
	\subfigure[Curve by RPIA(10)]{
		\includegraphics[width=0.4\textwidth]{Ex3_Curve_m20000_n500-eps-converted-to.pdf}}
\caption{The cubic B-spline fitting curves given by RPIA with $m=20000$ and $n= 500$ for Example \ref{ex:RPIA3}.}
\label{fig:ex3-RPIA1+Result}
\end{figure}

\begin{figure}[!htb]
\centering
 	\subfigure[Curve by RPIA(5)]{
		\includegraphics[width=0.4\textwidth]{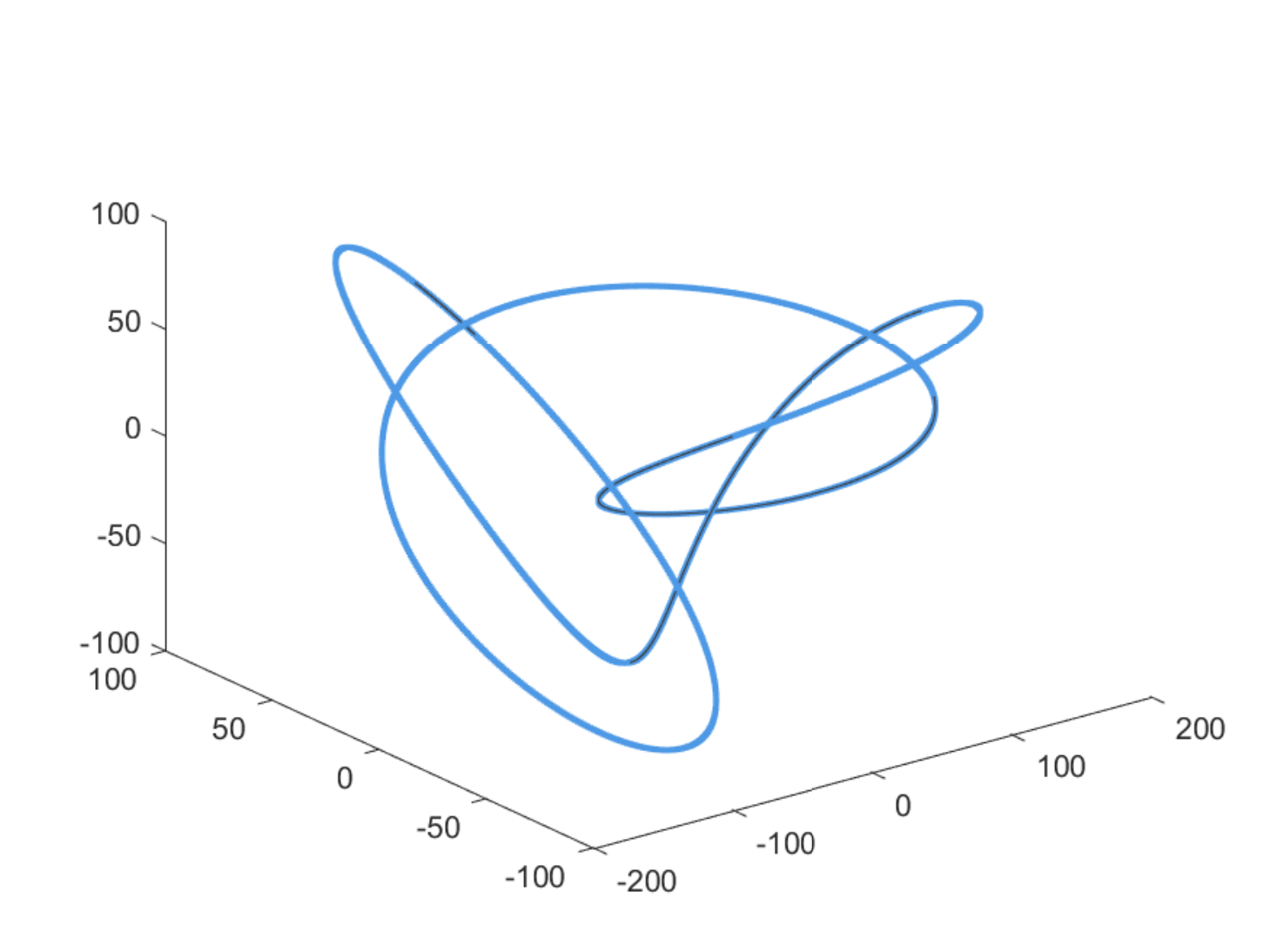}}
	\subfigure[Curve by RPIA(10)]{
		\includegraphics[width=0.4\textwidth]{Ex4_Curve_m20000_n500-eps-converted-to.pdf}}
\caption{The cubic B-spline fitting curves given by RPIA with $m=20000$ and $n= 500$ for Example \ref{ex:RPIA4}.}
\label{fig:ex4-RPIA1+Result}
\end{figure}

\begin{figure}[!htb]
\centering
	\subfigure[Example \ref{ex:RPIA1}]{
		\includegraphics[width=0.4\textwidth]{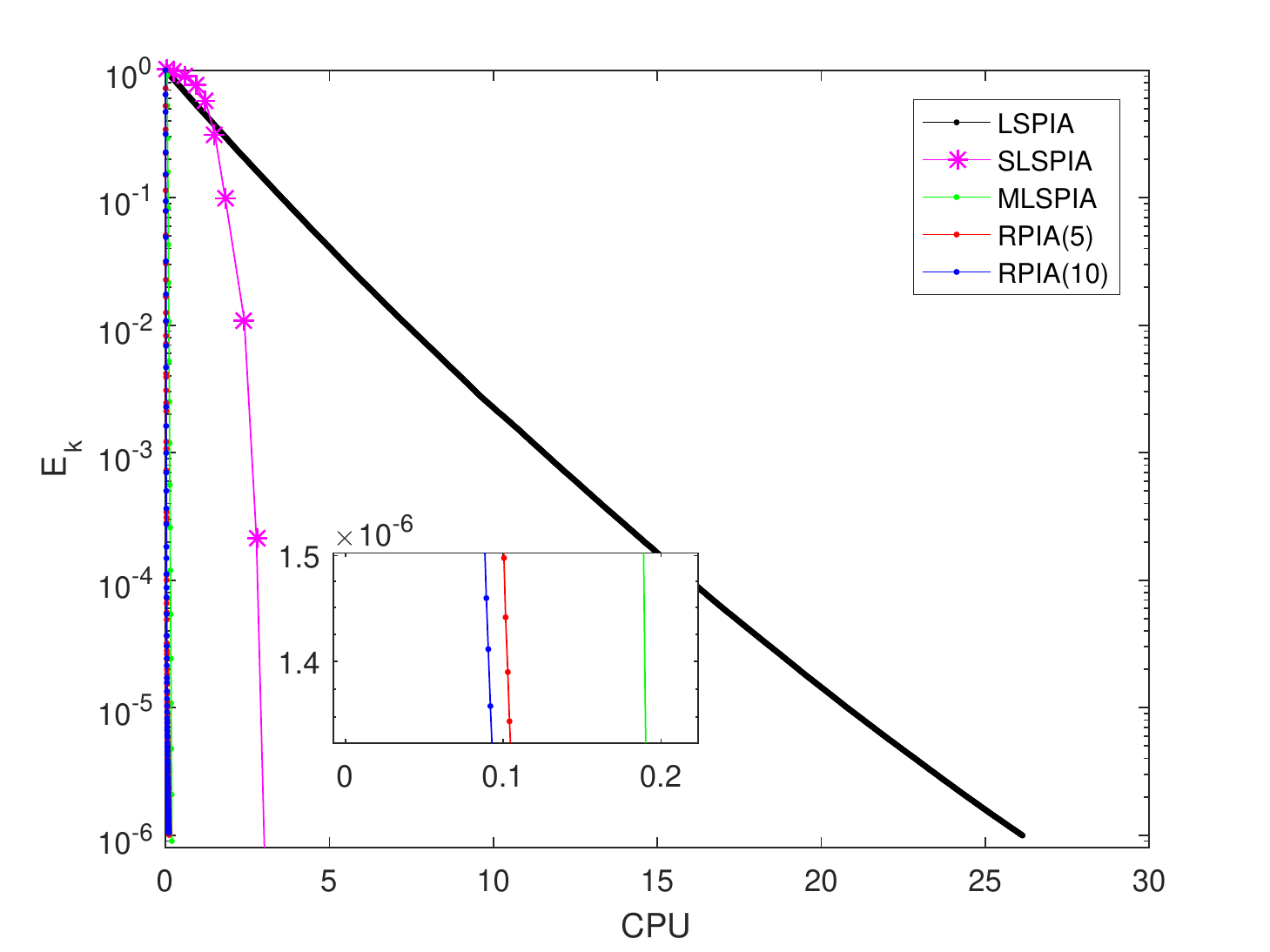}}
	\subfigure[Example \ref{ex:RPIA2}]{
		\includegraphics[width=0.4\textwidth]{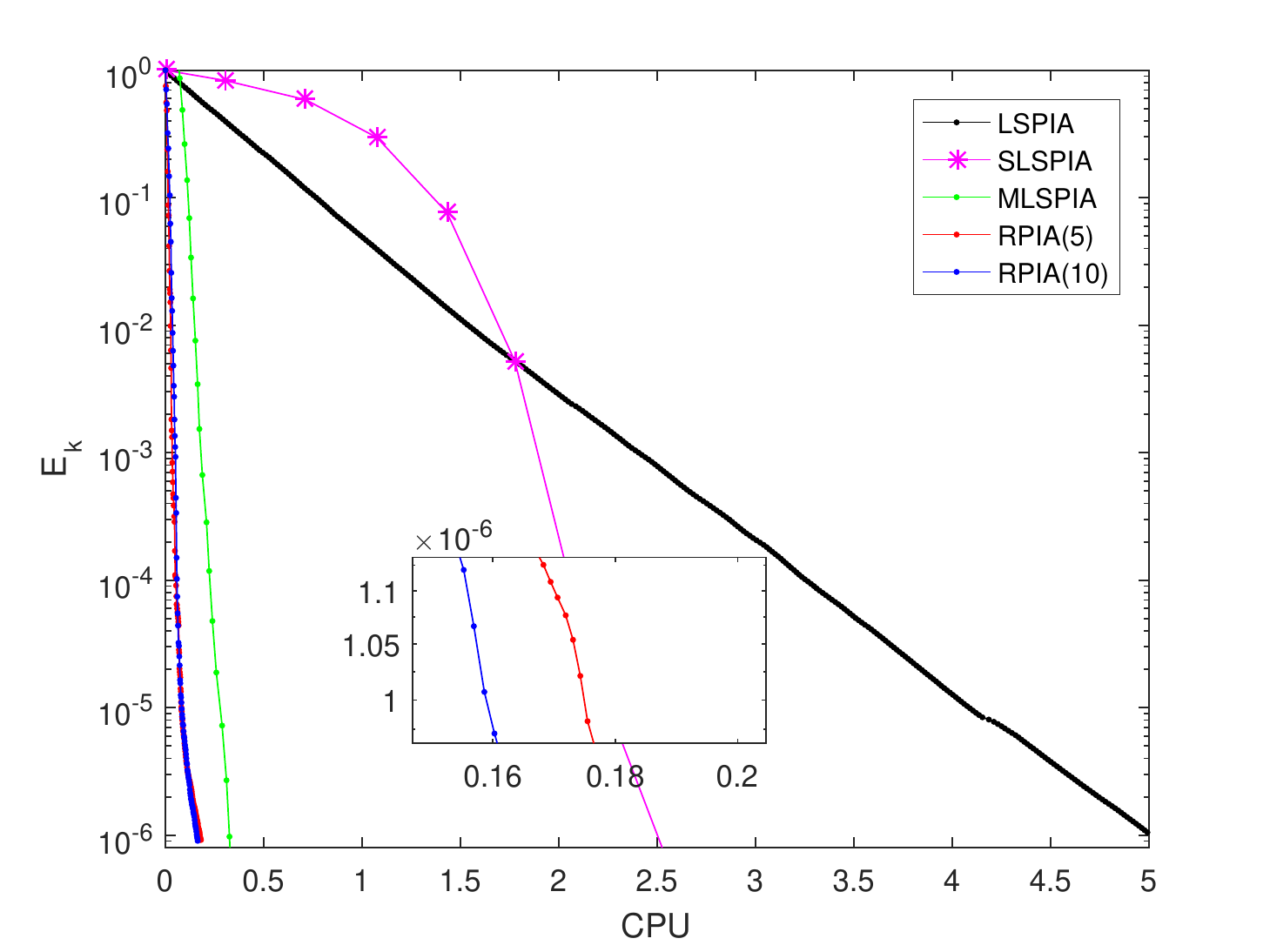}}
	\subfigure[Example \ref{ex:RPIA3}]{
        \includegraphics[width=0.4\textwidth]{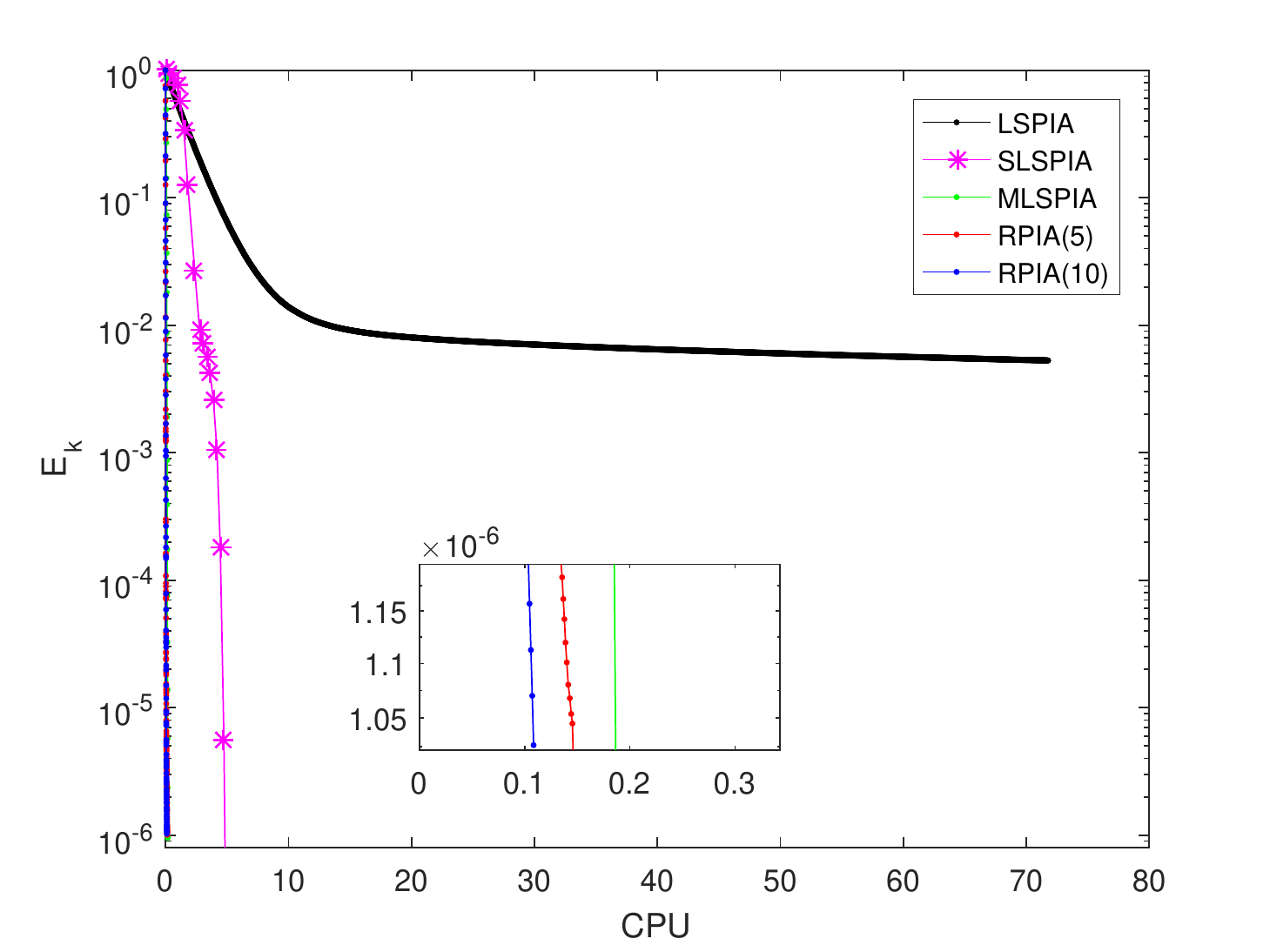}}
	\subfigure[Example \ref{ex:RPIA4}]{
        \includegraphics[width=0.4\textwidth]{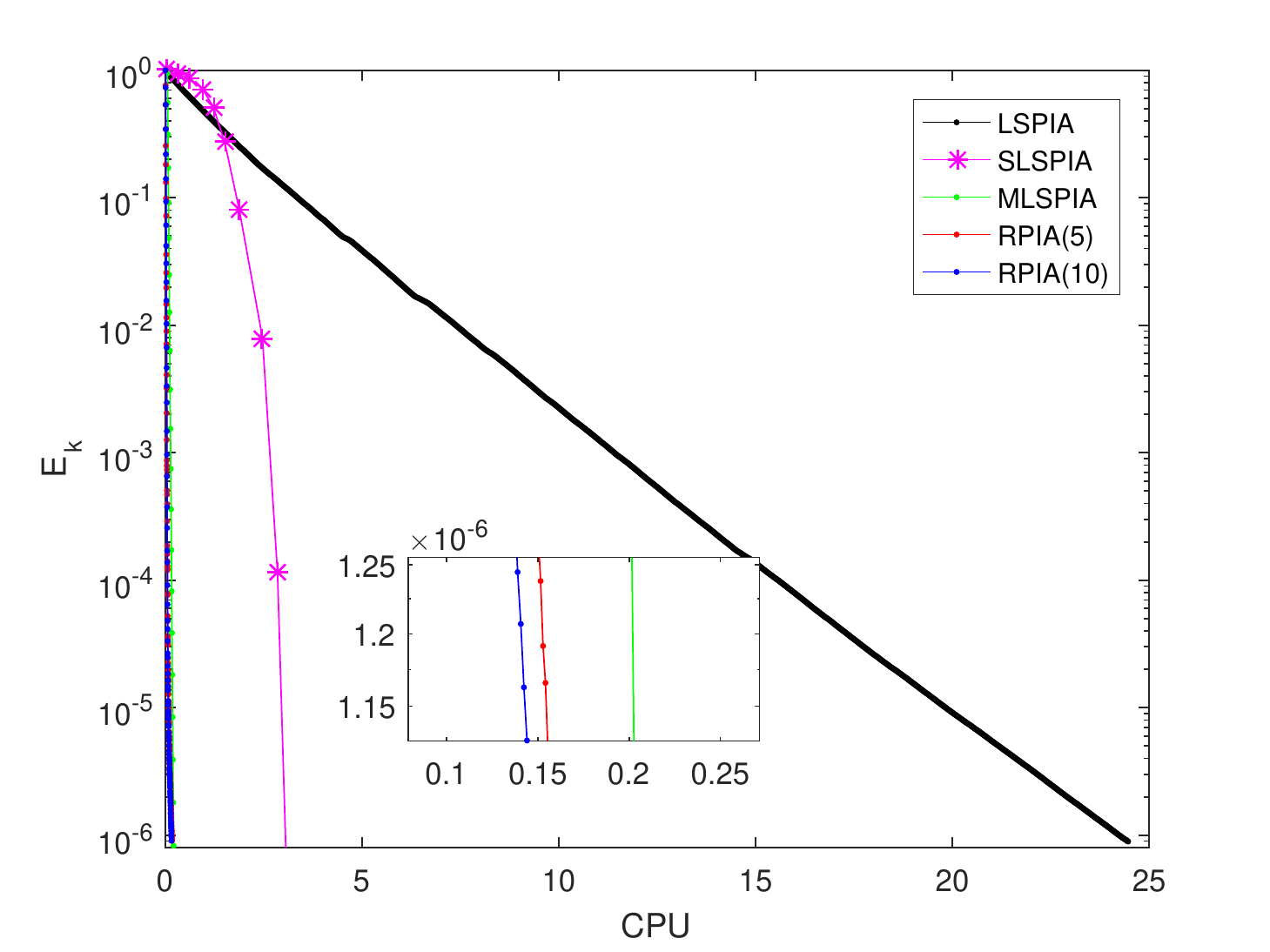}}
\caption{The iteration histories of $E_k$  versus  CPU for LSPIA, SLSPIA, MLSPIA, RPIA(5), and RPIA(10) using cubic B-spline curves with $m=20000$ and $n= 500$  for Examples \ref{ex:RPIA1} (a), \ref{ex:RPIA2} (b), \ref{ex:RPIA2} (c), and \ref{ex:RPIA4} (d).}
\label{fig:ex1-ex4+FittingError}
\end{figure}

The numerical results are listed in Tables  \ref{tab:ex1-RPIA1+Result}--\ref{tab:ex4-RPIA1+Result}. We find that the  relative fitting errors of the four methods are comparable while RPIA takes  much less time than LSPIA, SLSPIA, and MLSPIA in all settings. In Figures \ref{fig:ex1-RPIA1+Result}--\ref{fig:ex4-RPIA1+Result}, we draw the curves constructed by RPIA when the numbers of data and control points are $20000$ and $500$, respectively. Accordingly, Figure \ref{fig:ex1-ex4+FittingError} shows the iteration history of relative fitting error for the tested methods. It is clear that the relative fitting errors of RPIA decay faster than that of MLSPIA and much faster than that of LSPIA and SLSPIA when the computing time increases, which indicates that RPIA is more effective than LSPIA, SLSPIA, and MLSPIA  in actual applications.

\begin{table}[!htb]
 \normalsize
\caption{Example  \ref{ex:RPIA1}: $E_\infty$, IT, and CPU for LSPIA, SLSPIA, MLSPIA, RPIA(5), and RPIA(10) with $n=500$ and various $m$.}
\centering
\begin{tabular}{ ccccccc}
\cline{1-7}
&&LSPIA&SLSPIA&MLSPIA&RPIA(5)&RPIA(10)\\
\cline{1-7}
$m=20000$&{$E_\infty$}&$8.92 \times 10^{-7}$&$3.38 \times 10^{-7}$&$3.86 \times 10^{-7}$&$8.90 \times 10^{-7}$&$8.90 \times 10^{-7}$\\
         &{IT}&1866&10&27&6118.2&3753.1\\
         &{CPU}&26.251&3.046&0.198&0.112&0.098\\
\cline{1-7}
$m=30000$&{$E_\infty$}&$8.89 \times 10^{-7}$&$3.39 \times 10^{-7}$&$8.74 \times 10^{-7}$&$8.90 \times 10^{-7}$&$8.90 \times 10^{-7}$\\
         &{IT}&1867&10&27&5781.8&3533.4\\
         &{CPU}&40.765&6.321&0.317&0.183&0.140\\
\cline{1-7}
\end{tabular}
\label{tab:ex1-RPIA1+Result}
\end{table}

\begin{table}[!htb]
 \normalsize
\caption{Example  \ref{ex:RPIA2}: $E_\infty$, IT, and CPU for LSPIA, SLSPIA, MLSPIA, RPIA(5), and RPIA(10) with $n=500$ and various $m$.}
\centering
\begin{tabular}{ ccccccc}
\cline{1-7}
&&LSPIA&SLSPIA&MLSPIA&RPIA(5)&RPIA(10)\\
\cline{1-7}
$m=20000$&{$E_\infty$}&$8.64 \times 10^{-7}$&$7.03 \times 10^{-7}$&$3.43 \times 10^{-7}$&$8.90 \times 10^{-7}$&$8.88 \times 10^{-7}$\\
         &{IT} & 336&8&25&6008.7&3849.5\\
         &{CPU}&5.061&2.786&0.305&0.153&0.138\\
\cline{1-7}
$m=30000$&{$E_\infty$}&$8.67 \times 10^{-7}$&$7.57 \times 10^{-7}$&$7.69 \times 10^{-7}$&$8.90 \times 10^{-7}$&$8.90 \times 10^{-7}$\\
         &{IT} & 336&8&25&6001.7&3807.4\\
         &{CPU}&7.817&4.700&0.339&0.203&0.168\\
\cline{1-7}
\end{tabular}
\label{tab:ex2-RPIA1+Result}
\end{table}

\begin{table}[!htb]
 \normalsize
\caption{Example  \ref{ex:RPIA3}: $E_\infty$, IT, and CPU for LSPIA, SLSPIA, MLSPIA, RPIA(5), and RPIA(10) with $n=500$ and various $m$.}
\centering
\begin{tabular}{ ccccccc}
\cline{1-7}
&&LSPIA&SLSPIA&MLSPIA&RPIA(5)&RPIA(10)\\
\cline{1-7}
$m=20000$&{$E_\infty$}&$\ddag$&$5.31 \times 10^{-7}$&$3.80 \times 10^{-7}$&$8.92 \times 10^{-7}$&$8.92 \times 10^{-7}$\\
         &{IT}&$>10^{5}$&17&26&6133.4&4002.8\\
         &{CPU}&$\ddag$&5.694&0.199&0.160&0.107\\
\cline{1-7}
$m=30000$&{$E_\infty$}&$\ddag$&$5.34 \times 10^{-7}$&$8.35 \times 10^{-7}$&$8.92 \times 10^{-7}$&$8.91 \times 10^{-7}$\\
         &{IT}&$>10^{5}$&17&26&5949.3&3941.4\\
         &{CPU}&$\ddag$&7.561&0.290&0.133&0.118\\
\cline{1-7}
\end{tabular}
\begin{tablenotes}
\footnotesize
\item[1.] {\it The item ' $>10^{5}$' represents that the number of iteration steps exceeds $10^{5}$. In this case, the corresponding relative fitting error and CPU time are expressed as '$\ddag$'.}
\end{tablenotes}
\label{tab:ex3-RPIA1+Result}
\end{table}

\begin{table}[!htb]
 \normalsize
\caption{Example  \ref{ex:RPIA4}: $E_\infty$, IT, and CPU for LSPIA, SLSPIA, MLSPIA, RPIA(5), and RPIA(10) with $n=500$ and various $m$.}
\centering
\begin{tabular}{ ccccccc}
\cline{1-7}
&&LSPIA&SLSPIA&MLSPIA&RPIA(5)&RPIA(10)\\
\cline{1-7}
$m=20000$&{$E_\infty$}&$ 8.92 \times 10^{-7}$&$ 5.24 \times 10^{-7}$&$ 8.29 \times 10^{-7}$&$ 8.91 \times 10^{-7}$&$ 8.91 \times 10^{-7}$\\
         &{IT}&1649&10&27&6079.2&4178.6\\
         &{CPU}&24.462&3.173&0.206&0.177&0.150\\
\cline{1-7}
$m=30000$&{$E_\infty$}&$ 8.92 \times 10^{-7}$&$ 5.24 \times 10^{-7}$&$ 8.29 \times 10^{-7}$&$ 8.91 \times 10^{-7}$&$ 8.91 \times 10^{-7}$\\
         &{IT}        &1649&10&27&6074.9&4169.5\\
         &{CPU}       &38.143&5.869&0.344&0.181&0.159\\
\cline{1-7}
\end{tabular}
\label{tab:ex4-RPIA1+Result}
\end{table}

\subsection{Surface fitting}
\begin{figure}[!htb]
\centering
	\subfigure[Example \ref{ex:RPIA_surface1}]{
	\includegraphics[width=0.4\textwidth]{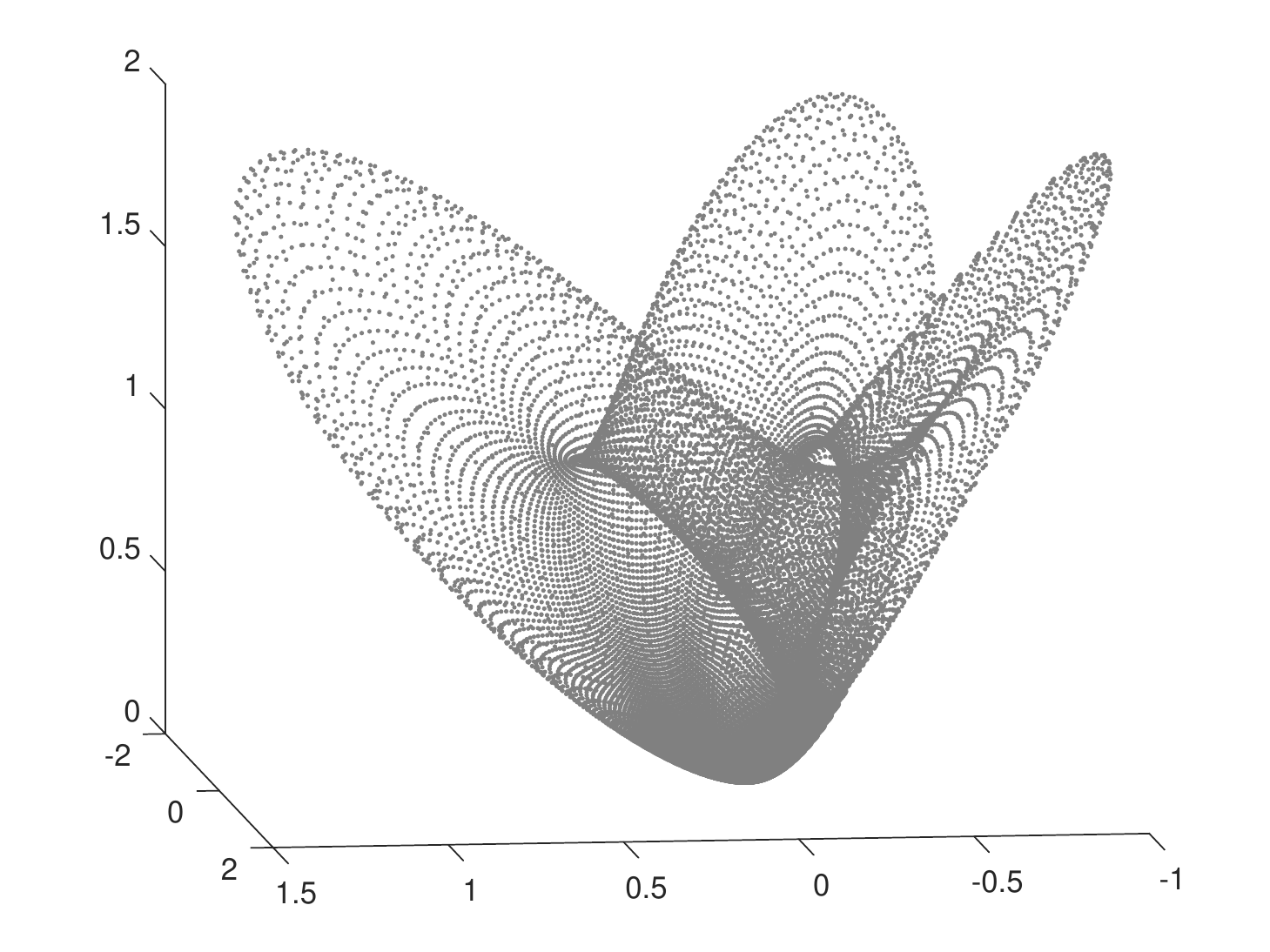}}
 	\subfigure[Example \ref{ex:RPIA_surface2}]{
	\includegraphics[width=0.4\textwidth]{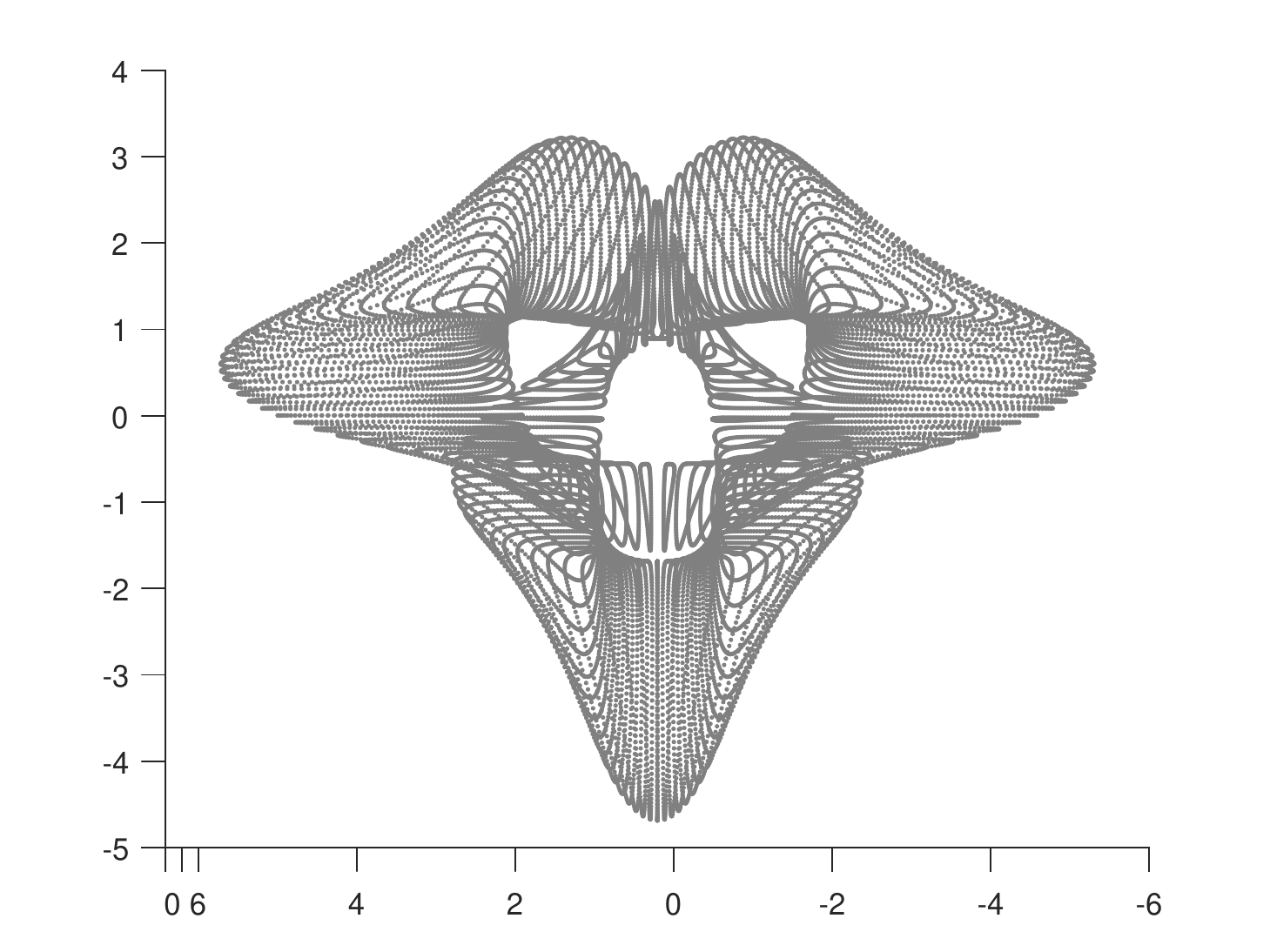}}
 	\subfigure[Example \ref{ex:RPIA_surface3}]{
	\includegraphics[width=0.4\textwidth]{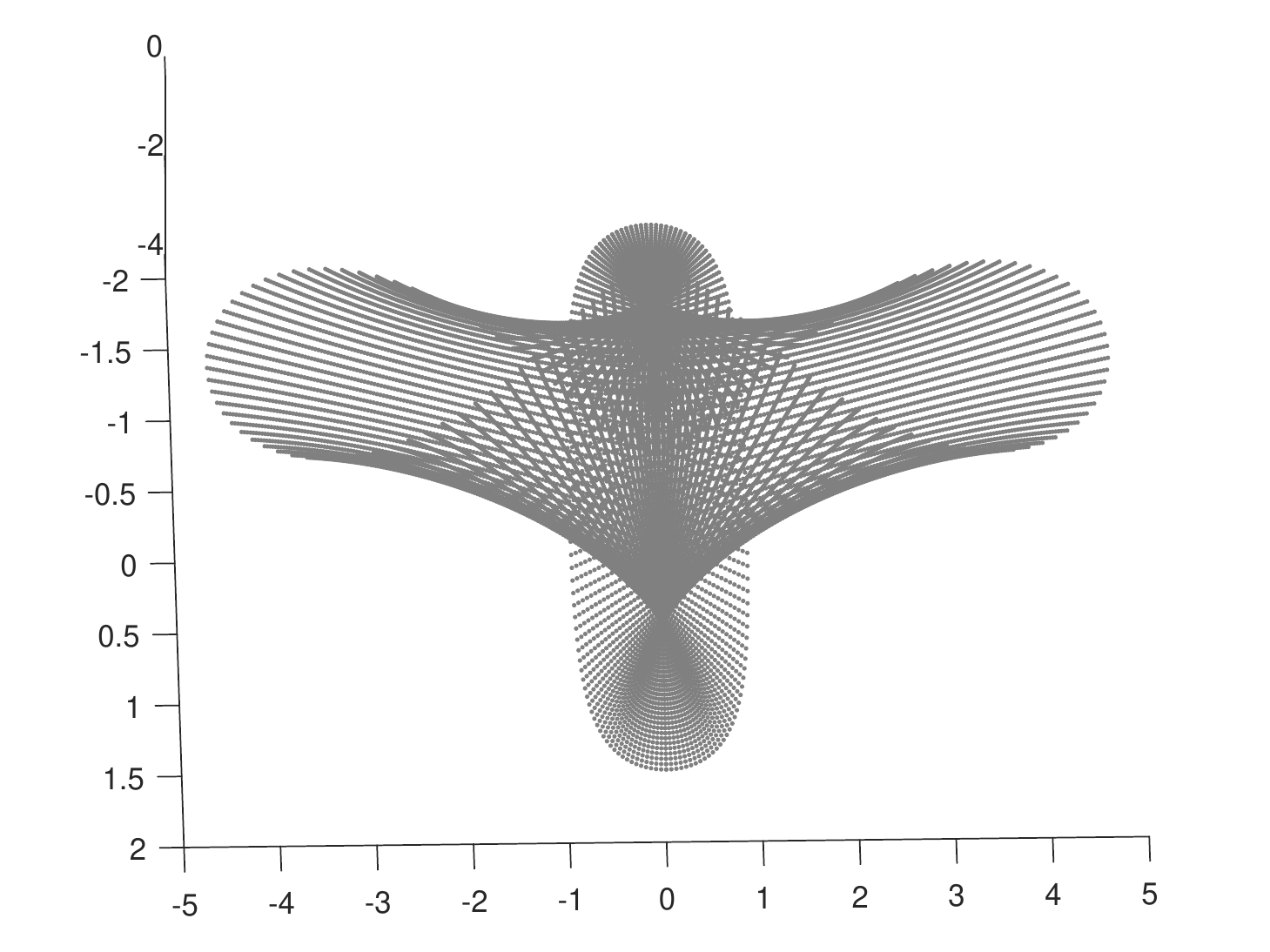}}
 	\subfigure[Example \ref{ex:RPIA_surface4}]{
	\includegraphics[width=0.4\textwidth]{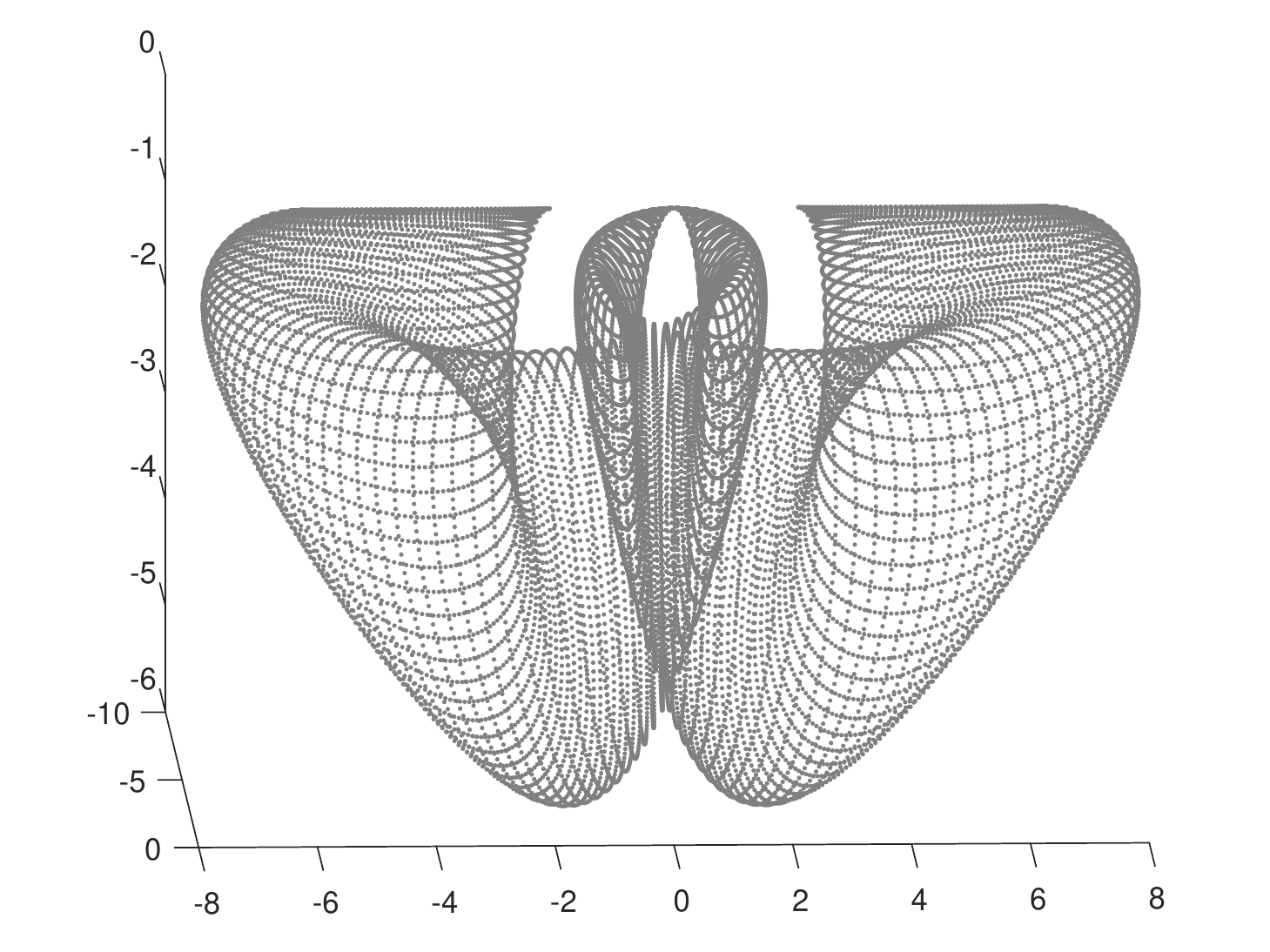}}
\caption{The data point sets to be fitted in Examples \ref{ex:RPIA_surface1} (a), \ref{ex:RPIA_surface2} (b), \ref{ex:RPIA_surface3} (c), and \ref{ex:RPIA_surface4} (d) with $m=p=120$.}
\label{fig:surface+initial}
\end{figure}

Similar to the case of curve fitting, we organize the execution details of RPIA for surface fitting as follows. We assign the parameters  $\left\{x_h \right\}_{h=0}^m$ and $\left\{y_l \right\}_{l=0}^p$ as
\begin{equation*}
 x_0 = 0,~x_h = x_{h-1}+\frac{\sum_{t=0}^p \bT{{\bm Q}_{ht}-{\bm Q}_{h-1, t}}}{\sum_{s=1}^m \sum_{t=0}^p \bT{{\bm Q}_{st}-{\bm Q}_{s-1, t}}},
~{\rm and}~
 x_m =1
\end{equation*}
for $h=1,2,\cdots,m$ and
\begin{equation*}
 y_0 = 0,~y_l = y_{l-1}+\frac{\sum_{s=0}^m \bT{{\bm Q}_{sl}-{\bm Q}_{s, l-1}}}{\sum_{h=0}^m \sum_{t=1}^p \bT{{\bm Q}_{st}-{\bm Q}_{s, t-1}}},
~{\rm and}~
 y_p =1
\end{equation*}
for $l=1,2,\cdots,p$, respectively, and define two knot vectors as
\begin{equation*}
  \left\{
  0,~0,~0,~0,
  ~\bar{x}_4,~\bar{x}_5,\cdots,\bar{x}_n,
  ~1,~1,~1,~1
  \right\},
\end{equation*}
  with $\bar{x}_{h+3}=(1-\alpha_x)x_{i-1} + \alpha_x x_i$, $i = \lfloor hd_x \rfloor$, $\alpha_x = hd_x-i$, $d_x = (m+1)/(n-2)$ for $h=1,2,\cdots,n-3$ and
\begin{equation*}
  \left\{
  0,~0,~0,~0,
  ~\bar{y}_4,~\bar{y}_5,\cdots,\bar{y}_n,
  ~1,~1,~1,~1
  \right\},
\end{equation*}
  with $\bar{y}_{l+3}=(1-\alpha_y)y_{j-1} + \alpha_y y_j$, $j = \lfloor ld_y \rfloor$, $\alpha_y=ld_y-j$, $d_y = (p+1)/(n-2)$ for $l=1,2,\cdots,n-3$, respectively. The initial control points are selected as
${\bm P}^{(0)}_{ij} =  {\bm Q}_{f_1(i), f_2(j)}$,
where $f_2(0)=0$, $f_2(n)=p$, $f_2(j)=\lfloor  pj/n  \rfloor$  for $j=1,2,\cdots,n-1$, and the definition of $f_1(i)$ is given by Section \ref{subsec:NE1}.  We consider the following four point sets, as shown in Figure \ref{fig:surface+initial}.

\vskip 1ex
\begin{example}\label{ex:RPIA_surface1}
$(m+1)\times (p+1)$ data points sampled uniformly from a boy surface, whose parametric equation is given by
\begin{equation*}
\left\{ \begin{array}{l}
\bx = 2/3 (\cos(t) \cos(2t) + \sqrt{2}\sin(t)\cos(s))\cos(t)/(\sqrt{2} - \sin(2t)\sin(3s)), \vspace{1ex}\\
\by = 2/3 (\cos(t) \sin(2t) - \sqrt{2}\sin(t)\sin(s))\cos(t)/(\sqrt{2} - \sin(2t)\sin(3s)),  \vspace{1ex}\\
\bz = \sqrt{2}\cos(t)\cos(t)/(\sqrt{2} - \sin(2t)\sin(3s))
~ (-\pi \leq t,~s \leq  \pi).
\end{array}\right .
\end{equation*}
\end{example}

\begin{example}\label{ex:RPIA_surface2}
$(m+1)\times (p+1)$ data points sampled uniformly from a tranguloid-trefoil surface, whose parametric equation is given by
\begin{equation*}
\left\{ \begin{array}{l}
\bx = 2 \sin(3t) /(2 + \cos(s)),                      \vspace{1ex}\\
\by = 2( \sin(t) + 2\sin(2t)) / (2 + \cos(s+2\pi/3)), \vspace{1ex}\\
\bz =  (\cos(t) - 2\cos(2t)) (2 + \cos(s)) (2 + \cos(s + 2\pi/3))/4
~ (-\pi \leq t,~s \leq  \pi).
\end{array}\right .
\end{equation*}
\end{example}

\begin{example}\label{ex:RPIA_surface3}
$(m+1)\times (p+1)$ data points sampled uniformly from a Verrill-minimal surface, whose parametric equation is given by
\begin{equation*}
\left\{ \begin{array}{l}
\bx = -2 t \cos(s) + 2 \cos(s) / t - 2 t^3 \cos(3 s) / 3,\vspace{1ex}\\
\by = 6t\sin(s) - 2 \sin(s) /t - 2 t^3 \sin(3 s) / 3,\vspace{1ex}\\
\bz = 4 \log(t)
~ (0.5 \leq t \leq 1,~0 \leq s \leq 2 \pi).
\end{array}\right .
\end{equation*}
\end{example}

\begin{example}\label{ex:RPIA_surface4}
$(m+1)\times (p+1)$ data points sampled uniformly from a bent-horns surface, whose parametric equation is given by
\begin{equation*}
\left\{ \begin{array}{l}
\bx = (2 + \cos(t)) (s/3 - \sin(s)),         \vspace{1ex}\\
\by = (2 + \cos(t - 2 \pi / 3)) (\cos(s) - 1), \vspace{1ex}\\
\bz = (2 + \cos(t + 2 \pi / 3)) (\cos(s) - 1)
~(-\pi \leq t \leq  \pi,~-2\pi \leq s \leq 2\pi).
\end{array}\right .
\end{equation*}
\end{example}

\begin{table}[!htb]
 \normalsize
\caption{Example  \ref{ex:RPIA_surface1}: $E_\infty$, IT, and CPU for LSPIA, SLSPIA, MLSPIA, RPIA(5), and RPIA(10) with $n=20$ and various $m$ and $p$.}
\centering
\begin{tabular}{ ccccccc}
\cline{1-7}
&&LSPIA&SLSPIA&MLSPIA&RPIA(5)&RPIA(10)\\
\cline{1-7}
$m=p=120$&{$E_\infty$}&$ 8.09\times 10^{-7}$&$ 5.03\times 10^{-7}$&$7.76 \times 10^{-7}$&$8.77 \times 10^{-7}$&$8.89 \times 10^{-7}$\\
         &{IT}& 28& 10& 124& 6952.7 & 5419.0 \\
         &{CPU}& 0.261& 1.662& 0.671& 0.139& 0.133\\
\cline{1-7}
$m=p=160$&{$E_\infty$}&$8.66 \times 10^{-7}$&$5.48 \times 10^{-7}$&$8.49 \times 10^{-7}$&$8.88 \times 10^{-7}$&$8.88 \times 10^{-7}$\\
         &{IT}& 28& 10& 127& 6708.2&5239.2 \\
         &{CPU}& 0.612&3.279 &1.336 & 0.143& 0.139 \\
\cline{1-7}
\end{tabular}
\label{tab:ex5-RPIA1+Result}
\end{table}

\begin{table}[!htb]
 \normalsize
\caption{Example  \ref{ex:RPIA_surface2}: $E_\infty$, IT, and CPU for LSPIA, SLSPIA, MLSPIA, RPIA(5), and RPIA(10) with $n=20$ and various $m$ and $p$.}
\centering
\begin{tabular}{ ccccccc}
\cline{1-7}
&&LSPIA&SLSPIA&MLSPIA&RPIA(5)&RPIA(10)\\
\cline{1-7}
$m=p=120$&{$E_\infty$}&$8.61 \times 10^{-7}$&$4.77 \times 10^{-7}$&$8.47 \times 10^{-7}$&$8.84 \times 10^{-7}$&$8.86 \times 10^{-7}$\\
         &{IT}& 37& 11& 169 & 6530.0&4979.7 \\
         &{CPU}& 0.375& 1.908&0.975 & 0.191&0.178 \\
\cline{1-7}
$m=p=160$&{$E_\infty$}&$8.45 \times 10^{-7}$&$4.23 \times 10^{-7}$&$8.59 \times 10^{-7}$&$8.85 \times 10^{-7}$&$8.81 \times 10^{-7}$\\
         &{IT}&38 &11 & 114& 6458.8& 5138.1\\
         &{CPU}&0.725 &3.239 &1.194 & 0.197&0.182 \\
\cline{1-7}
\end{tabular}
\label{tab:ex6-RPIA1+Result}
\end{table}

\begin{table}[!htb]
 \normalsize
\caption{Example  \ref{ex:RPIA_surface3}: $E_\infty$, IT, and CPU for LSPIA, SLSPIA, MLSPIA, RPIA(5), and RPIA(10) with $n=20$ and various $m$ and $p$.}
\centering
\begin{tabular}{ ccccccc}
\cline{1-7}
&&LSPIA&SLSPIA&MLSPIA&RPIA(5)&RPIA(10)\\
\cline{1-7}
$m=p=120$&{$E_\infty$}&$8.90 \times 10^{-7}$&$3.37 \times 10^{-7}$&$8.41 \times 10^{-7}$&$8.86 \times 10^{-7}$&$8.85 \times 10^{-7}$\\
         &{IT}&45 &12 &227 &6829.4 &5220.7 \\
         &{CPU}&0.483 &2.173 & 1.335&0.148 &0.132 \\
\cline{1-7}
$m=p=160$&{$E_\infty$}&$8.60 \times 10^{-7}$&$8.05 \times 10^{-7}$&$8.10 \times 10^{-7}$&$8.86 \times 10^{-7}$&$8.85 \times 10^{-7}$\\
         &{IT}&21 & 9& 137&6559.3 &5027.5 \\
         &{CPU}&0.395 &2.773 &1.432 &0.145 &0.140 \\
\cline{1-7}
\end{tabular}
\label{tab:ex7-RPIA1+Result}
\end{table}

\begin{table}[!htb]
 \normalsize
\caption{Example  \ref{ex:RPIA_surface4}: $E_\infty$, IT, and CPU for LSPIA, SLSPIA, MLSPIA, RPIA(5), and RPIA(10) with $n=20$ and various $m$ and $p$.}
\centering
\begin{tabular}{ ccccccc}
\cline{1-7}
&&LSPIA&SLSPIA&MLSPIA&RPIA(5)&RPIA(10)\\
\cline{1-7}
$m=p=120$&{$E_\infty$}&$8.45 \times 10^{-7}$&$6.44 \times 10^{-7}$&$8.04 \times 10^{-7}$&$8.84 \times 10^{-7}$&$8.87 \times 10^{-7}$\\
         &{IT}&34 &11 &143 &6803.7 &4861.8 \\
         &{CPU}&0.328 &2.016 &0.826 &0.176 &0.168 \\
\cline{1-7}
$m=p=160$&{$E_\infty$}&$7.61 \times 10^{-7}$&$4.34 \times 10^{-7}$&$8.14 \times 10^{-7}$&$8.81 \times 10^{-7}$&$8.89 \times 10^{-7}$\\
         &{IT}&20 & 9&106 &6618.5 &4932.5 \\
         &{CPU}&0.387 &3.066 &1.224 &0.207 &0.178 \\
\cline{1-7}
\end{tabular}
\label{tab:ex8-RPIA1+Result}
\end{table}

\begin{figure}[!htb]
\centering
 	\subfigure[Surface by RPIA(5)]{
		\includegraphics[width=0.4\textwidth]{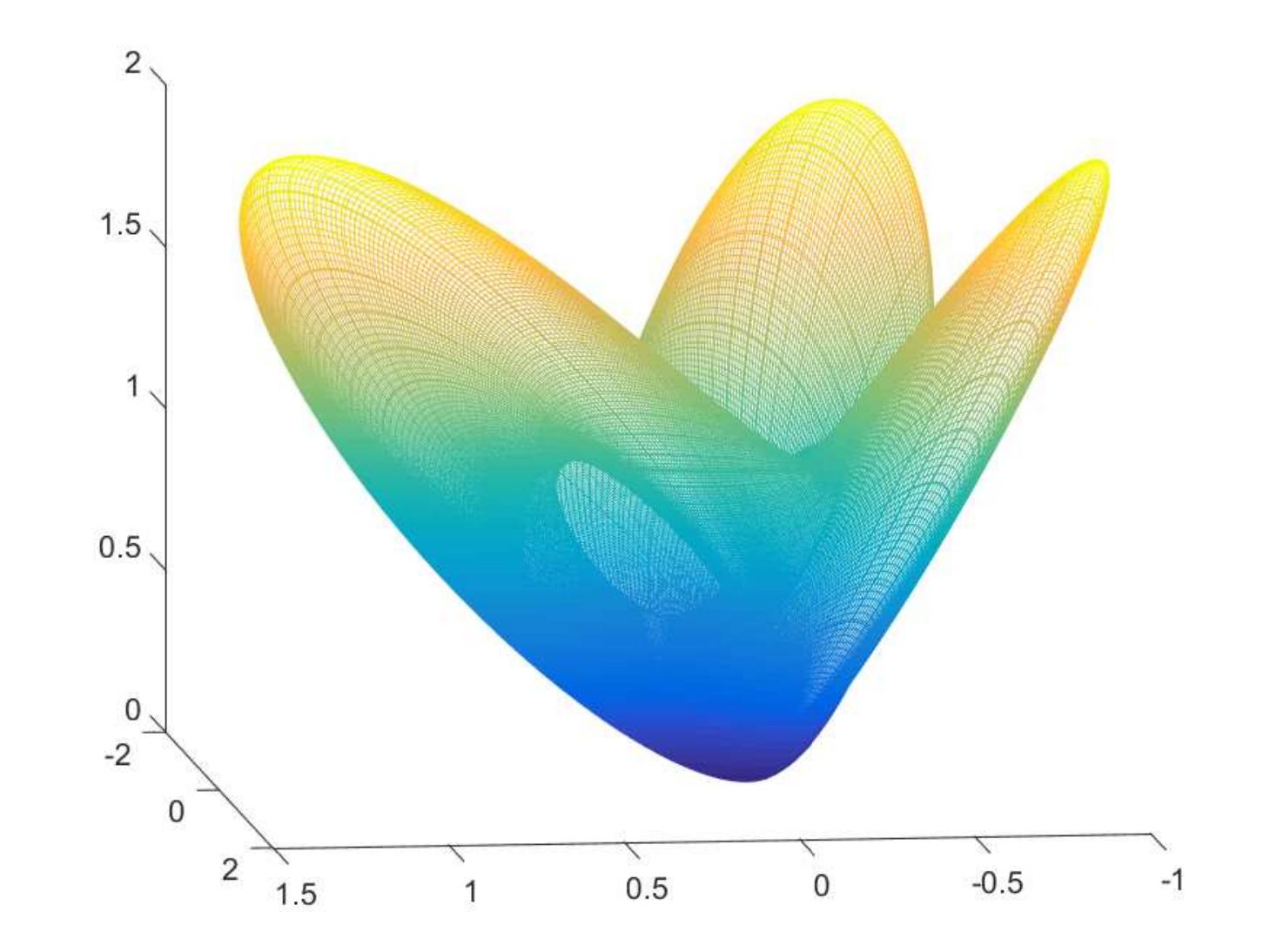}}\vspace{10pt}
	\subfigure[Surface by RPIA(10)]{
	    \includegraphics[width=0.4\textwidth]{Surface_m150_p150_Ex5-eps-converted-to.pdf}}
\caption{The bi-cubic B-spline fitting surfaces given by RPIA with $m=p=120$ and $n=20$ for Example \ref{ex:RPIA_surface1}.}
\label{fig:ex5-RPIA1+Result}
\end{figure}

\begin{figure}[!htb]
\centering
 	\subfigure[Surface by RPIA(5)]{
		\includegraphics[width=0.4\textwidth]{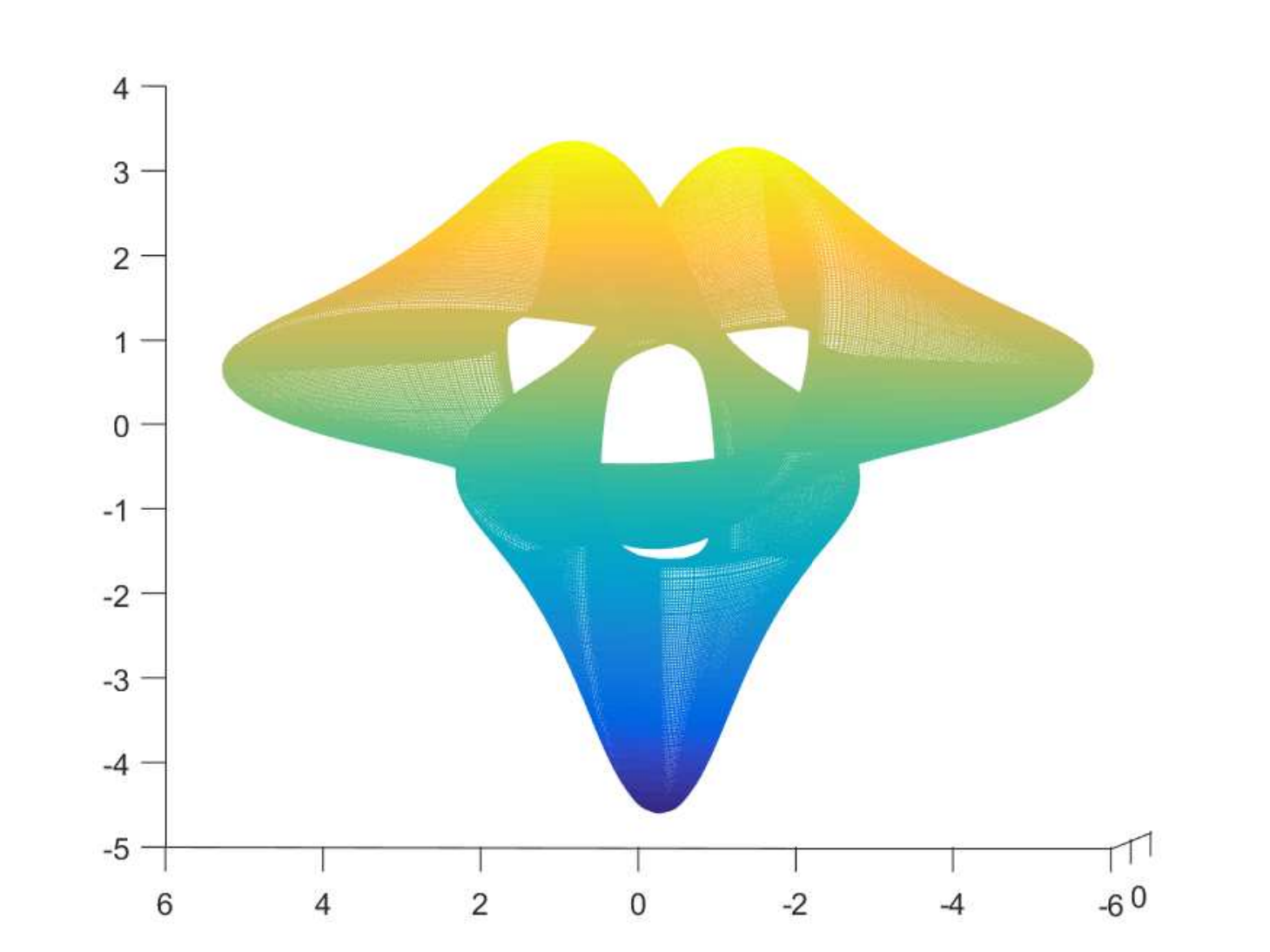}}\vspace{10pt}
	\subfigure[Surface by RPIA(10)]{
	    \includegraphics[width=0.4\textwidth]{Surface_m150_p150_Ex6-eps-converted-to.pdf}}
\caption{The bi-cubic B-spline fitting surfaces given by RPIA with $m=p=120$ and $n=20$  for Example \ref{ex:RPIA_surface2}.}
\label{fig:ex6-RPIA1+Result}
\end{figure}

\begin{figure}[!htb]
\centering
 	\subfigure[Surface by RPIA(5)]{
		\includegraphics[width=0.4\textwidth]{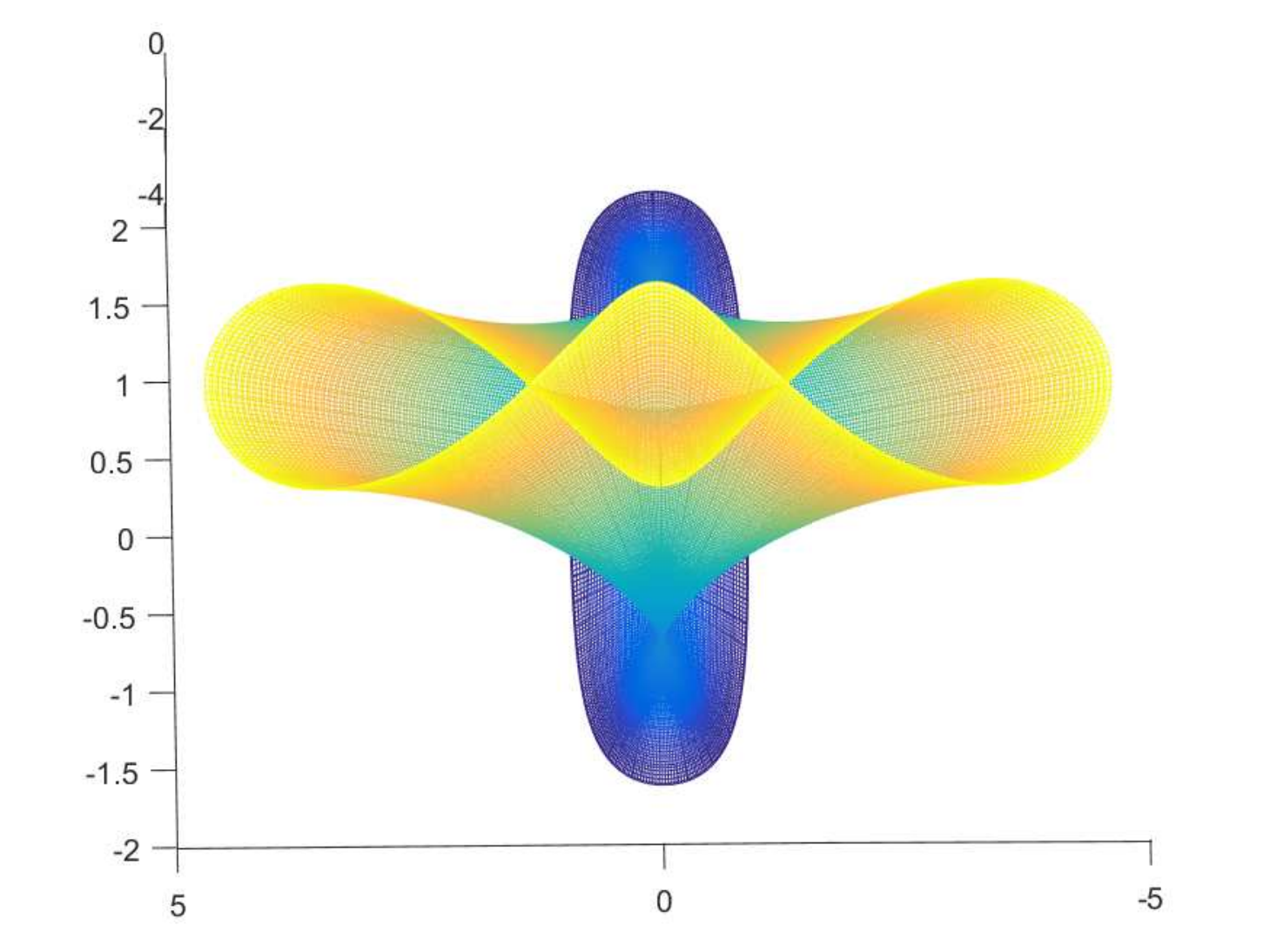}}
	\subfigure[Surface by RPIA(10)]{
	    \includegraphics[width=0.4\textwidth]{Surface_m150_p150_Ex7-eps-converted-to.pdf}}
\caption{The bi-cubic B-spline fitting surfaces given by RPIA with $m=p=120$ and $n=20$  for Example \ref{ex:RPIA_surface3}.}
\label{fig:ex7-RPIA1+Result}
\end{figure}

\begin{figure}[!!htb]
\centering
 	\subfigure[Surface by RPIA(5)]{
		\includegraphics[width=0.4\textwidth]{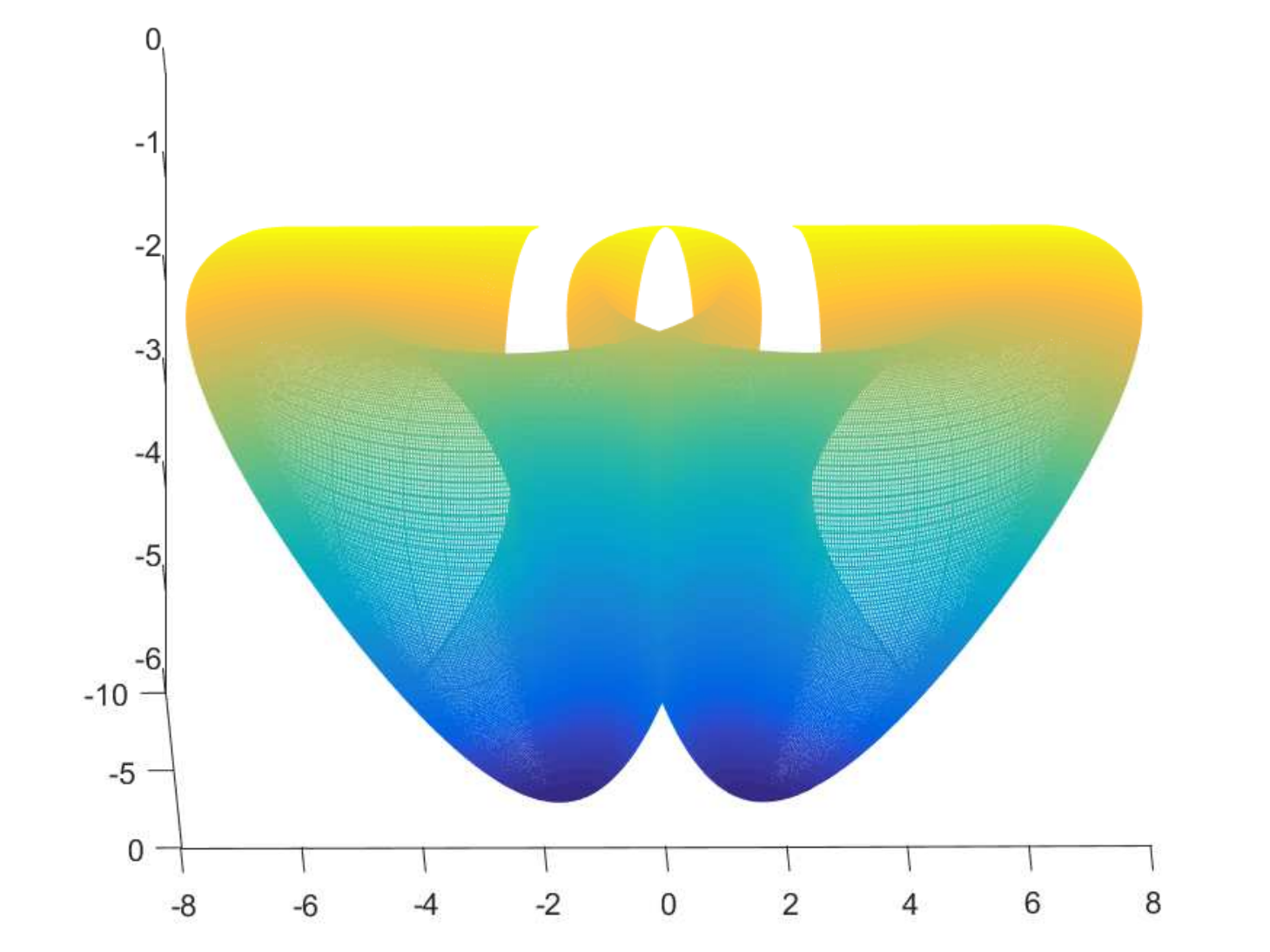}}
	\subfigure[Surface by RPIA(10)]{
	    \includegraphics[width=0.4\textwidth]{Surface_m150_p150_Ex8-eps-converted-to.pdf}}
\caption{The bi-cubic B-spline fitting surfaces given by RPIA with $m=p=120$ and $n=20$  for Example \ref{ex:RPIA_surface4}.}
\label{fig:ex8-RPIA1+Result}
\end{figure}

\begin{figure}[!htb]
\centering
	\subfigure[Example \ref{ex:RPIA_surface1}]{ \includegraphics[width=0.4\textwidth]{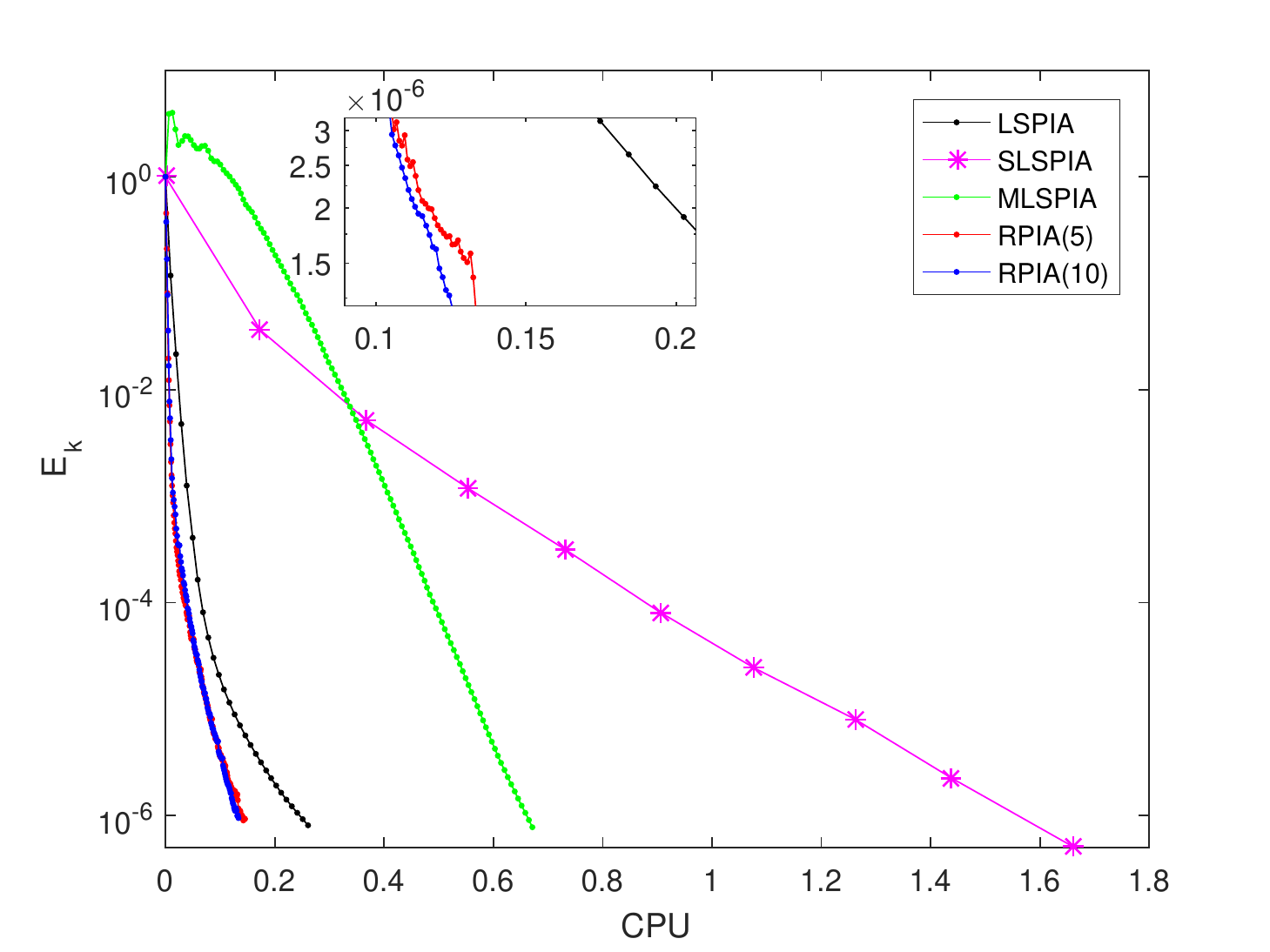}}
	\subfigure[Example \ref{ex:RPIA_surface2}]{ \includegraphics[width=0.4\textwidth]{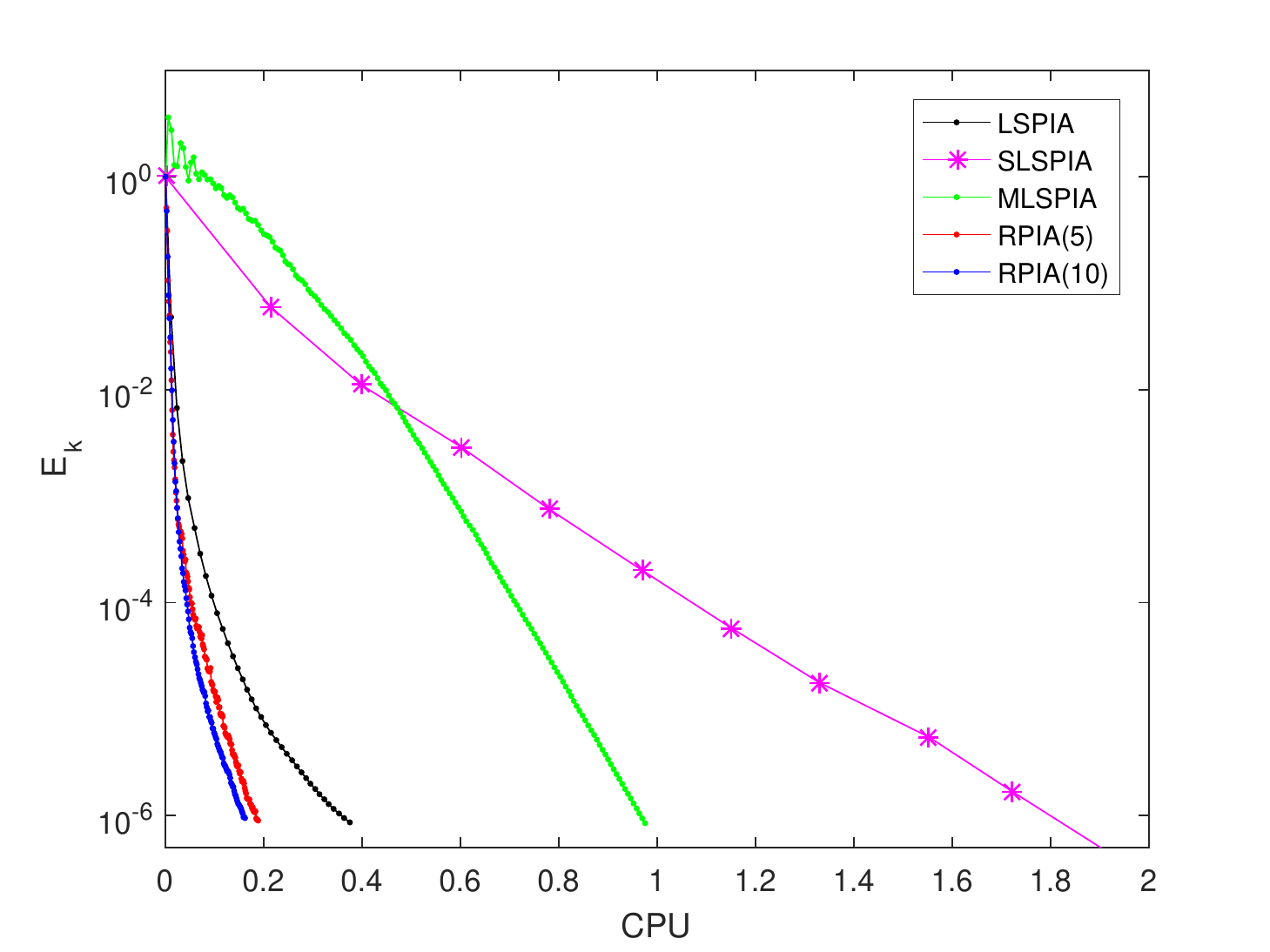}}
	\subfigure[Example \ref{ex:RPIA_surface3}]{ \includegraphics[width=0.4\textwidth]{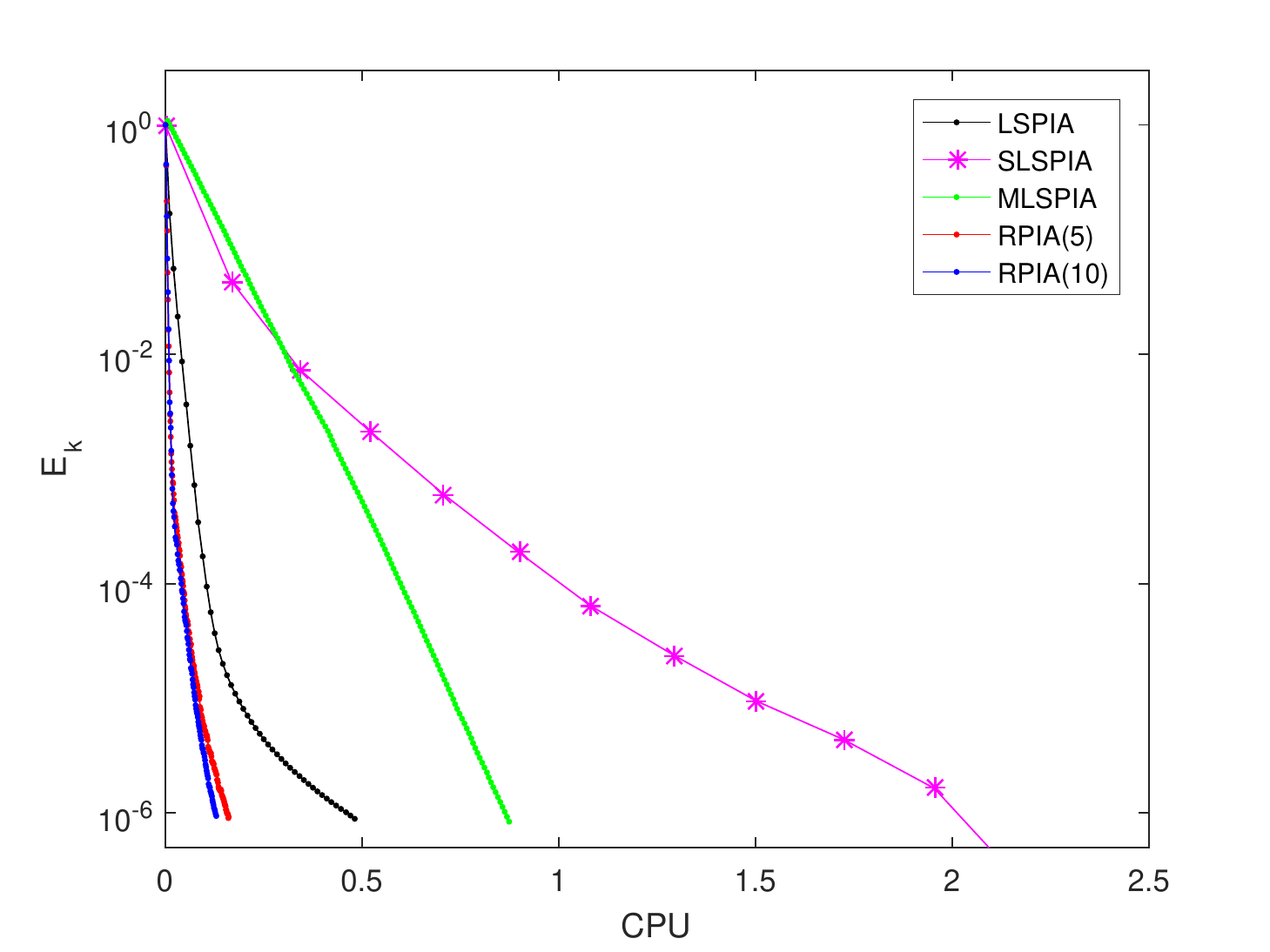}}
	\subfigure[Example \ref{ex:RPIA_surface4}]{ \includegraphics[width=0.4\textwidth]{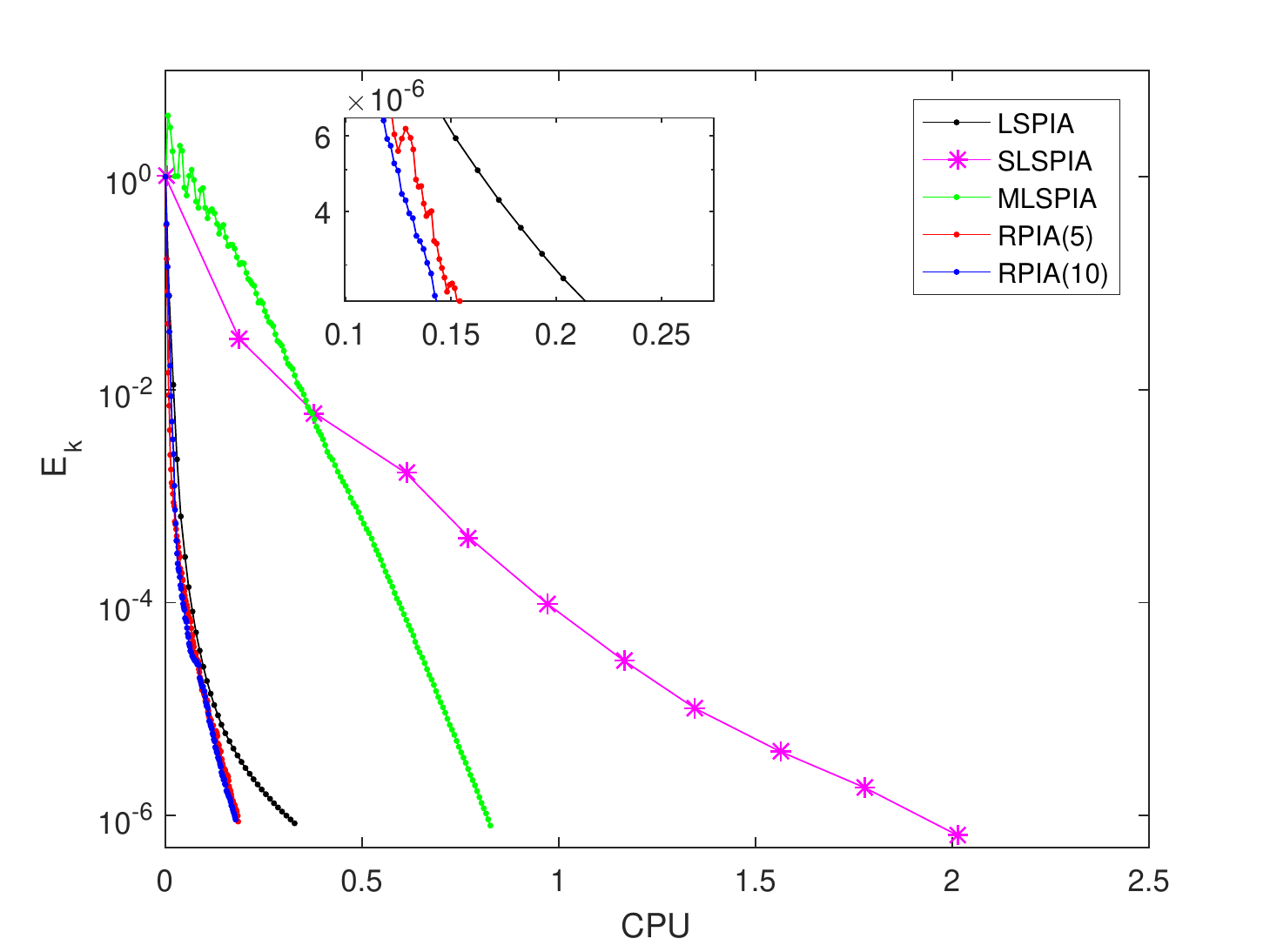}}
\caption{The iteration histories of $E_k$  versus  CPU for LSPIA, SLSPIA, MLSPIA, RPIA(5), and RPIA(10) using bi-cubic B-spline surfaces with $m=p=120$ and $n=20$  for Examples \ref{ex:RPIA_surface1} (a), \ref{ex:RPIA_surface2} (b), \ref{ex:RPIA_surface3} (c), and \ref{ex:RPIA_surface4} (d).}
\label{fig:ex5-ex8+FittingError}
\end{figure}

Tables \ref{tab:ex5-RPIA1+Result}--\ref{tab:ex8-RPIA1+Result} show the numerical results of $E_\infty$, IT, and CPU for LSPIA, SLSPIA, MLSPIA, RPIA(5), and RPIA(10) with $n=20$ and different $m$ and $p$ for Examples \ref{ex:RPIA_surface1}--\ref{ex:RPIA_surface4}. It can be seen that the IT for RPIAs with the block size $5$ and $10$ are more than that of LSPIA, SLSPIA, and MLSPIA, but the CPU for RPIAs are less than that of these methods. It indicates that RPIA is more efficient than LSPIA, SLSPIA, and MLSPIA. Figures \ref{fig:ex5-RPIA1+Result}--\ref{fig:ex8-RPIA1+Result} show the surfaces constructed by RPIA when  $m=p=150$ and $n=20$. At the same time,  we also draw the iteration history of relative fitting error for the tested methods in Figure \ref{fig:ex5-ex8+FittingError}. It is shown that the relative fitting errors of RPIA decay much faster than that of LSPIA, SLSPIA, and MLSPIA when the computing time increases.

\section{Concluding remarks}\label{Con}\noindent
The column form of LSPIA implies that it will operate on all of the control points at each iteration. In this work, inspired by this discovery, we turn to update the partial control points according to a randomized index set and keep the other ones remaining unchanged, and propose the randomized progressive iterative approximation to fit data. Our approach includes linear algebra-based procedures, the principle of PIA, and the promising randomized index selection method. The advantages of using the randomized technique include that: the resulting algorithm is easier to analyze and implement, has lower memory requirements, and is more parallelizable in practice.  We prove that our method has local PIA property and obtains a least-squares result in the limit sense.  From algebraic aspects, for RPIA curve and surface fittings, it is equivalent to use the randomized block coordinate descent method to solve the linear systems of the form ${\bm M}{\bm x}={\bm g}$ and ${\bm M}{\bm X}{\bm N} ={\bm G}$, respectively, which to the best of our knowledge have not been previously studied. We give some numerical examples to demonstrate the convergence behaviors of such a randomized iterative method. The numerical results show that RPIA is more effective than some often used LSPIA-type variants.

In the end some comments on the related recent works and possible extensions of RPIA are made as follows.

\begin{enumerate}[(1)]
\setlength{\itemindent}{0cm}
\vskip 0.5ex
\item We would like to mention a very important by-product from our convergence analysis, namely, it suggests a weighted RPIA iterative algorithm, like MLSPIA \cite{20HW} for LSPIA \cite{14DL}. In curve fitting, for example, by introducing a weight parameter in formula \eqref{Curve+eq:Update+deltak}, the control point is updated by
    \begin{align*}
    {\bm p}^{(k+1)}_{i_k} = {\bm p}^{(k)}_{i_k} +  \omega_{k}
    \frac{\sum_{j=0}^m \mu_{i_k}(x_j) {\bm r}_{j}^{(k)}}{ \sum_{i_k \in \mathrm{I}_{\imath_k}} \sum_{j=0}^{m}\mu_{i_k}^2(x_j)}
    \end{align*}
    for $k=0,1,2,\cdots$. The choice of an effective weight parameter is essential to accelerate the convergence of RPIA. In particular, if $\omega_{k}$ is one,  it recovers the standard RPIA method.
\vskip 0.5ex
\item The subdivision surface problems may arise from very different applications and modeling, e.g., computer graphics and feature film industry. Subdivision surface refers to a class of modeling schemes that define an object through recursive subdivision starting from an initial control mesh. Another feature of this work is that extend the RPIA method to approximate the vertices of a mesh by using Loop \cite{87Loop} and Catmull-Clark \cite{78CC} surfaces.
\end{enumerate}

\section*{Acknowledgment}
This work is supported by the National Natural Science Foundation of China under grants 12101225 and 12201651.

\end{document}